\documentclass[preprint,english]{siamltex}
\usepackage[T1]{fontenc}
\usepackage[latin9]{inputenc}
\usepackage{float}
\usepackage{amsmath}
\usepackage{amssymb}
\usepackage{graphicx}
\usepackage{xcolor}

\makeatletter

\newcommand{\lyxdot}{.}
\floatstyle{ruled}
\newfloat{algorithm}{tbp}{loa}
\providecommand{\algorithmname}{Algorithm}
\floatname{algorithm}{\protect\algorithmname}

\usepackage{fullpage}

\usepackage{algorithm,algpseudocode}

\newcommand{\change}[1]{#1}

\makeatother

\usepackage{babel}

\def\la{\langle}
\def\ra{\rangle}

\title{A low-rank projector-splitting integrator \\for the Vlasov--Poisson equation}

\author{Lukas Einkemmer\footnotemark[1]\ \footnotemark[3] \and Christian Lubich\footnotemark[1]}

\begin{document}

\maketitle

\renewcommand{\thefootnote}{\fnsymbol{footnote}}
\footnotetext[1]{Mathematisches Institut,
       Universit\"at T\"ubingen,
       Auf der Morgenstelle 10,
       D--72076 T\"ubingen,
       Germany. Email: {\tt \{einkemmer,lubich\}@na.uni-tuebingen.de}}
 \footnotetext[3]{Department of Mathematics, University of Innsbruck, Austria}
\renewcommand{\thefootnote}{\arabic{footnote}}

\begin{abstract} Many problems encountered in plasma physics require a description by kinetic equations, which are posed in an up to six-dimensional phase space. A direct discretization of this phase space, often called the Eulerian approach, has many advantages but is extremely expensive from a computational point of view.
In the present paper we propose a dynamical low-rank approximation to the Vlasov--Poisson equation, with time integration by a particular splitting method. This approximation is derived by constraining the dynamics to a manifold of low-rank functions via a tangent space projection and by splitting this projection into the subprojections from which it is built. This reduces a time step for the six- (or four-) dimensional Vlasov--Poisson equation to  solving two systems of three- (or two-) dimensional advection equations over the time step, once in the position variables and once in the velocity variables, where the size of each system of advection equations is equal to the chosen rank. By a hierarchical dynamical low-rank approximation, a time step for the Vlasov--Poisson equation can be further reduced to a set of six (or four) systems of one-dimensional advection equations, where the size of each system of advection equations is still equal to the rank. The resulting systems of advection equations can then be solved by standard techniques such as semi-Lagrangian or spectral methods. Numerical simulations in two and four dimensions for linear Landau damping, for a two-stream instability \change{and for a plasma echo problem} highlight the favorable behavior of this numerical method and show that the proposed algorithm is able to drastically reduce the required computational effort.
\end{abstract}

\begin{keywords} Dynamical low-rank approximation, projector-splitting integrator, Vlasov--Poisson equation.
\end{keywords}

\section{Introduction}

Many physical phenomena in both space
and laboratory plasmas require a kinetic description. The typical
example of such a kinetic model is the Vlasov\textendash Poisson equation\begin{equation}\label{eq:vlasov-poisson}
\begin{aligned}
& \partial_{t}f(t,x,v)+v\cdot\nabla_{x}f(t,x,v)-E(f)(x)\cdot\nabla_{v}f(t,x,v)=0 \\
& \nabla\cdot E(f)(x)=-\int f(t,x,v)\,\mathrm{d}v+1,\qquad\;\;\nabla\times E(f)(x)=0,
\end{aligned}
\end{equation}
\change{which models the time evolution of electrons in a collisionless plasma in the electrostatic
regime (assuming a constant background ion density).} Equation \eqref{eq:vlasov-poisson} has to be supplemented
with appropriate boundary and initial conditions. The particle-density
function $f(t,x,v)$, where $x\in\Omega_{x}\subset\mathbb{R}^{d}$ and
$v\in\Omega_{v}\subset\mathbb{R}^{d}$ (with spatial dimension $d\le 3$), 
is the quantity to be computed. 

In most applications, the primary computational challenge stems from
the fact that the equation is posed as an evolution equation in an up to six-dimensional phase
space. Thus, direct discretization of the Vlasov\textendash Poisson equation on
a grid requires the storage of $\mathcal{O}(n^{2d})$ floating point numbers, where $n$
is the number of grid points per direction.
This has been, and in many cases still is, prohibitive
from a computational point of view. 

Consequently, particle methods are widely
used in many application areas that require the numerical solution
of kinetic models. Particle methods, such as the particle-in-cell
approach, only discretize the physical space $x$ (but not the
velocity space $v$). The $v$-dependence is approximated by initializing
a large number of particles that follow the characteristic curves
of equation \eqref{eq:vlasov-poisson}. This can significantly reduce
the computational cost and has the added benefit that particles congregate
in high-density regions of the phase space. However, particle methods
suffer from numerical noise that only decreases as the square root
in the number of particles. This is, in particular, an issue for problems
where the tail of the distribution function has to be resolved accurately.
For a detailed discussion the reader is referred to the review
article \cite{verboncoeur2005particle}.

As an alternative, the entire $2d$-dimensional phase space can be
discretized. This Eulerian approach has received significant
attention in recent years. Due to the substantial increase in computational
power, four-dimensional simulation can now be performed
routinely on high performance computing (HPC) systems. In addition,
some five and six-dimensional simulations (usually where in one direction
a small number of grid points is sufficient) have been carried out.
Nevertheless, these simulations are still extremely costly from a
computational point of view. Consequently a significant research
effort has been dedicated to improving both the numerical algorithms
used (see, for example, \cite{crouseilles2011,einkemmer2014,einkemmer2016,filbet2003,grandgirard,klimas1994,Morrison2017,qiu2011,rossmanith2011,sonnendrucker1999,cheng1976,crouseilles2017,sircombe2009valis,crouseilles2015hamiltonian,crouseilles2016asymptotic})
as well as the corresponding parallelization to state of the art HPC
systems (see, for example, \cite{rozar2013,einkemmer2015,bigot2013,latu2007gyrokinetic,mehrenberger2013vlasov,einkemmer2016mixed,crouseilles2009parallel}).

The seminal paper  \cite{cheng1976} introduced a
time-splitting approach. This approach has the advantage that the
nonlinear Vlasov\textendash Poisson equation can be reduced to a sequence
of one-dimensional advection equations \change{(see, for example, \cite{fijalkow1999,filbet2001,einkemmer2018comparison})}. Splitting schemes that have
similar properties have been proposed for a range of more complicated models
(including the relativistic Vlasov\textendash Maxwell equations \change{\cite{sircombe2009valis,crouseilles2015hamiltonian}}, drift-kinetic models \change{\cite{grandgirard,crouseilles2017}}, and gyrokinetic models \change{\cite{grandgirard2016}}). This is then usually combined with
a semi-Lagrangian discretization of the resulting advection equations.
A semi-Lagrangian approach has the advantage that the resulting numerical
method is completely free from a Courant\textendash Friedrichs\textendash Lewy
(CFL) condition. 

Nevertheless, all these methods still incur a computational and storage
cost that scales as $\mathcal{O}(n^{2d})$. Thus, it seems natural
to ask whether any of the dimension reduction techniques developed
for high-dimensional problems are able to help alleviate the computational burden of
the Eulerian approach. To understand why it is not entirely clear
if such an approach could succeed, we will consider the following
equation (with $x,v\in \mathbb{R}$)
\[
\partial_{t}f(t,x,v)+v\,\partial_{x}f(t,x,v)=0.
\]
This is just the one-dimensional Vlasov\textendash Poisson equation
\eqref{eq:vlasov-poisson} without the nonlinear term. If the initial
value $f_{0}(x,v)=\cos(kx)\mathrm{e}^{-v^{2}/2}$ with $k\in\mathbb{R}$
is imposed, then the exact solution can be easily written down as
\begin{equation}
f(t,x,v)=\cos(k(x-vt))\,\mathrm{e}^{-v^{2}/2}.\label{eq:toy-solution}
\end{equation}
From this expression we see that, as the system evolves in time, the
wave number in the $v$-direction is given by $kt$. In particular,
we have $\Vert\partial_{v}^{m}f(t,\cdot,\cdot))\Vert\propto(kt)^{m}$
which implies that, even though the solution is infinitely often differentiable,
smaller and smaller scales appear in phase space. This phenomenon
is referred to as filamentation and in many problems a reasonable
resolution of these small scale structures is important to obtain
physically meaningful results.

A specific problem for sparse-grid approximations is that mixed derivatives
(between $x$ and $v$, but perhaps even more problematic between
different $v$-directions in a multidimensional setting) can be huge.
Nevertheless, sparse-grid techniques have been considered for the
Vlasov\textendash Poisson equation \cite{Kormann,Deriaz2016,Guo2016}.
A further challenge for these methods is that evaluating the particle-density
function (which is required for semi-Lagrangian methods) can be expensive. 

Despite these drawbacks, sparse grids can be advantageous in some
situations. For example, in \cite{Kormann} first a tensor product
decomposition in the $x$- and $v$-directions is applied. The remaining
$d$-dimensional problems are then approximated using a sparse grid.
For Landau damping in four dimensions this scheme was able to reduce
the required memory by approximately a factor of $10$ and the run
time by a factor of $2$. On the other hand, it is somewhat worrying
that, as was shown in the cited work, instabilities can
develop if not enough degrees of freedom are used. Furthermore, it
is well known that sparse grids do not resolve Gaussians very well.
This is a problem for kinetic simulations as the steady state is usually
a linear combination of Gaussians. In the previously mentioned work,
this deficiency is remedied by explicitly incorporating a Gaussian
into the numerical method (a so-called $\delta f$ scheme). However,
this approach does not work very well in cases where the particle density
significantly deviates from this form.

We now turn to the alternative approach taken in this paper, which appears to be novel for kinetic equations.
Computing numerical solutions of high-dimensional evolutionary partial differential
equations by  {\it dynamical low-rank approximation} has, however, been investigated extensively in quantum dynamics; see, in particular, 
 \cite{meyer90tmc,meyer09mqd} for the MCTDH approach to molecular quantum dynamics in the chemical physics literature and \cite{Lubich2008,lubich15tii,Conte2010} for a mathematical point of view to this approach. Some uses of dynamical low-rank approximation in areas outside quantum mechanics are described in \cite{Nonnenmacher2008,jahnke2008dynamical,Mena2017,Musharbash2018}. In a general mathematical setting, dynamical low-rank approximation has been studied in \cite{Koch2007,Koch2010,lubich13dab,arnold2014approximation}. A major algorithmic advance for the time integration was achieved with the projector-splitting methods first proposed in \cite{Lubich2014} for matrix differential equations and then developed further for various tensor formats in \cite{lubich15tii,lubich15tio,Haegeman2016,Kieri2016,Lubich2017}. In the present paper we will adapt such methods to the Vlasov\textendash Poisson equation.
 
For the Vlasov\textendash Poisson equation, low-rank
approximation to the solution is an interesting
option for a number of reasons. First, the solution of our toy problem,
given in equation (\ref{eq:toy-solution}), can be readily written
as
\[
f(t,x,v)=\left(\cos(kx)\cos(kvt)+\sin(kx)\sin(kvt)\right)\mathrm{e}^{-v^{2}/2}.
\]
Thus, the linear part of the equation can be represented exactly with
rank two. While this property is lost once we consider
the full Vlasov\textendash Poisson equation, it still gives some indication
that at least certain types of filamentation can be handled efficiently
by the numerical method based on a low-rank approximation. Second,
the typical kinetic equilibria formed by linear combination of Gaussians
can be very efficiently represented by a low-rank structure. This
is in stark contrast to the sparse-grid approach (as has been discussed
above). 

Recently, in  \cite{Kormanna} a low-rank numerical algorithm for
the Vlasov\textendash Poisson equation has been suggested that first
discretizes the problem in time and space. This is done using the
common technique of splitting in time and then applying a spline based
semi-Lagrangian scheme to the resulting advection equations. Each
step of such an algorithm can then be written as a linear combination
of low-rank approximations. To avoid that the rank grows during the
numerical simulation a singular value decomposition (SVD) is then
applied to truncate the low-rank approximation (i.e.,~approximate
it by a different low-rank representation with a fixed rank). 

In the present work we consider a different approach to approximating
the solution of the Vlasov\textendash Poisson equation by a low-rank
representation. We constrain the dynamics of the Vlasov\textendash Poisson
equation to a manifold of low-rank functions by a tangent space projection which is  then
split into its summands over a time step, adapting the projector-splitting approach to
time integration of \cite{Lubich2014}.
This  yields a sequence of advection equations in a lower-dimensional space. Then an appropriate
semi-Lagrangian scheme or even Fourier-based techniques can be applied
easily to obtain a numerical solution. This approach is first used to obtain
equations separately in $x$ and in $v$, which reduces a $2d$-dimensional problem
to a sequence of $d$-dimensional problems (this is the content of
Section \ref{sec:x-v decomposition}). If $r$ denotes the chosen rank, then over a time step one solves, one after the other,
\begin{itemize}
\item a system of $r$  advection equations in $x$ (and computing $2dr^2+r$ integrals over $\Omega_v$), 
\item a system of $r^2$ ordinary differential equations (and computing $2dr^2$ integrals over $\Omega_x$), 
\item a system of $r$ advection equations in $v$.
\end{itemize}
This first-order scheme can be symmetrized to yield a second-order time-stepping scheme. If $n$ degrees of freedom are used in each coordinate direction in the full discretization, then the storage cost is reduced to $O(rn^d)$ from $O(n^{2d})$, and similarly for the number of arithmetic operations.

If appropriate, the procedure
can then be repeated in a hierarchical manner in order to reduce the computations further to {\it one-dimensional}
advection equations and {\it one-dimensional} integrals (this is the content of Section \ref{sec:The-hierarchical-algorithm}).
We have chosen a hierarchical Tucker format since the dimensionality at
this point in the algorithm is only $d$ and this allows us to treat all coordinate directions in $x$ and in $v$ equally. The algorithm adapts the projector-splitting integrator for low-rank Tucker tensors from \cite{Lubich2017}.

A noteworthy feature of the proposed numerical scheme is that it
works entirely within the low-rank manifold.
No high-rank tensor needs to be formed and subsequently truncated.
This helps to reduce storage and computational cost. In addition,
depending on the problem either the physical space or the velocity
space or both can be solved directly, while still maintaining the
hierarchical low-rank approximation in the other parts of the algorithm.

Numerical simulations for linear Landau damping, a two-stream instability \change{and a plasma echo problem}
are presented in Section \ref{sec:Numerical-results}. There we show
that it is possible to use a small rank to capture the physics
of these problems. The proposed numerical method thus substantially
reduces the computational cost required to perform the simulation.
 
\section{Description of the numerical algorithm\label{sec:x-v decomposition}}

The goal of this section is to derive an algorithm that approximates
the solution of the Vlasov\textendash Poisson equation \eqref{eq:vlasov-poisson}
by a low-rank representation. To that end, the approximation to the particle-density function
$f(t,x,v)$ is constrained to the following form:
\begin{equation} 
f(t,x,v)\approx \sum_{i,j=1}^r X_{i}(t,x)S_{ij}(t)V_{j}(t,v),\label{eq:f-lowrank}
\end{equation}
where $S_{ij}(t)\in\mathbb{R}$ and we call $r$ the rank of the
representation. Note that the dependence of $f$ on the phase space
$(x,v)\in\Omega\subset\mathbb{R}^{2d}$ is now approximated by the
functions $\{X_{i}\colon i=1,\dots,r\}$ and $\{V_{j}\colon j=1,\dots,r\}$
which depend only on $x\in\Omega_{x}\subset\mathbb{R}^{d}$ and $v\in\Omega_{v}\subset\mathbb{R}^{d}$,
respectively, and on the time $t$. In the following
discussion we always assume that summation indices run from $1$ to
$r$ and we thus do not specify these bounds.

Now, we seek an approximation to the particle-density function
that for all $t$ lies in the set
\[
\overline{\mathcal{M}}=\biggl\{ f\in L^{2}(\Omega)\colon f(x,v)=\sum_{i,j}X_{i}(x)S_{ij}V_{j}(v)\text{ with }S_{ij}\in\mathbb{R},\,X_{i}\in L^{2}(\Omega_{x}),\,V_{j}\in L^{2}(\Omega_{v})\biggr\}.
\]
It is clear that this representation is not unique. In particular,
we can make the assumption that $\langle X_{i},X_{k}\rangle_x=\delta_{ik}$ and $\langle V_{j},V_{l}\rangle_v=\delta_{jl}$,
where $\la\cdot,\cdot\ra_x$ and $\la\cdot,\cdot\ra_v$ are the inner products on $L^{2}(\Omega_{x})$
and $L^{2}(\Omega_{v})$, respectively. We consider a path $f(t)$
on $\overline{\mathcal{M}}$. The corresponding time derivative is denoted
by $\dot{f}$ and is of the form 
\begin{equation}
\dot{f}=\sum_{i,j}\left(X_{i}\dot{S}_{ij}V_{j}+\dot{X}_{i}S_{ij}V_{j}+X_{i}S_{ij}\dot{V}_{j}\right).\label{eq:path}
\end{equation}
If we impose the gauge conditions $\la X_{i},\dot{X}_{j}\ra_x=0$ and $\la V_{i},\dot{V}_{j}\ra_v=0$,
then $\dot S_{ij}$ is uniquely determined by $\dot{f}$. This follows
easily from the fact that
\begin{equation}
\dot{S}_{ij}=\la X_{i}V_{j},\dot{f}\ra_{x,v}.\label{eq:unique-S}
\end{equation}
We then project both sides of equation (\ref{eq:path}) onto $X_{i}$
and $V_{j}$, respectively, and obtain
\begin{align}
\sum_{j}S_{ij}\dot{V}_{j} & =\la X_{i},\dot{f}\ra_x-\sum_{j}\dot{S}_{ij}V_{j},\label{eq:unique-V}\\
\sum_{i}S_{ij}\dot{X}_{i} & =\la V_{j},\dot{f}\ra_v-\sum_{i}X_{i}\dot{S}_{ij}.\label{eq:unique-X}
\end{align}
From these relations it follows that the $X_{i}$ and $V_{j}$ are
uniquely defined if the matrix $S=(S_{ij})$ is invertible. Thus, we seek an approximation
that for each time $t$ lies in the manifold
\begin{align*}
\mathcal{M} & =\biggl\{ f\in L^{2}(\Omega)\colon f(x,v)=\sum_{i,j}X_{i}(x)S_{ij}V_{j}(v)\text{ with invertible }S=(S_{ij})\in\mathbb{R}^{r\times r}, \\[-1mm]
&\qquad\qquad X_{i}\in L^{2}(\Omega_{x}),\,V_{j}\in L^{2}(\Omega_{v})
 \text{ with } \la X_{i},X_{k}\ra_x=\delta_{ik},\,\la V_{j},V_{l}\ra_v=\delta_{jl}\biggr\}
\end{align*}
with the corresponding tangent space 
\begin{align*}
\mathcal{T}_{f}\mathcal{M} & =\biggl\{\dot{f}\in L^{2}(\Omega)\colon\dot{f}(x,v)=\sum_{i,j}\left(X_{i}(x)\dot{S}_{ij}V_{j}(v)+\dot{X}_{i}(x)S_{ij}V_{j}(v)+X_{i}(x)S_{ij}\dot{V}_{j}(v)\right),\\
 & \qquad\qquad\text{with }\dot{S}\in\mathbb{R}^{r\times r},\,\dot{X}_{i}\in L^{2}(\Omega_{x}),\,\dot{V}_{j}\in L^{2}(\Omega_{v}),\text{ and }\la X_{i},\dot{X}_{j}\ra_x=0,\ \la V_{i},\dot{V_{j}}\ra_v=0\biggr\},
\end{align*}
where $f$ is given by equation (\ref{eq:f-lowrank}). Now, we consider
the reduced dynamics of the Vlasov\textendash Poisson equation on the manifold
$\mathcal{M}$. That is, we consider
for the approximate particle density (again denoted by $f$ for ease of notation)
\begin{equation}
\partial_{t}f=-P(f)\left(v\cdot\nabla_{x}f-E(f)\cdot\nabla_{v}f\right),\label{eq:vlasov-proj}
\end{equation}
where $P(f)$ is the orthogonal projector onto the tangent space $\mathcal{T}_{f}\mathcal{M}$,
as defined above.

We will consider the projection $P(f)g$ for a moment. From equations
(\ref{eq:path})-(\ref{eq:unique-X}) we obtain
\[
P(f)g=\sum_{j} \la V_{j},g\ra_v V_{j}-\sum_{i,j}X_{i} \la X_{i}V_{j},g\ra_{x,v} V_{j}+\sum_{i}X_i \la X_{i},g\ra_x.
\]
Let us introduce the two vector spaces $\overline{X}=\text{span}\left\{ X_{i}\colon i=1,\dots,r\right\} $
and $\overline{V}=\text{span}\left\{ V_{j}\colon j=1,\dots r\right\} $.
Then we can write the projector as follows

\begin{equation}
P(f)g=P_{\overline{V}}g-P_{\overline{V}}P_{\overline{X}}g+P_{\overline{X}}g,\label{eq:projector}
\end{equation}
where $P_{\overline X}$ is the orthogonal projector onto the vector space $\overline X$, which acts only on the $x$-variable of $g=g(x,v)$, and $P_{\overline V}$ is the orthogonal projector onto the vector space $\overline V$, which acts only on the $v$-variable of $g(x,v)$.
The decomposition of the projector into these three terms forms the
basis of our splitting procedure (for matrix differential equations this
has been first proposed in \cite{Lubich2014}). 

We proceed by substituting $g=v\cdot\nabla_{x}f-E(f)\cdot\nabla_{v}f$
into equation (\ref{eq:projector}). This immediately leads to a three-term
splitting for equation (\ref{eq:vlasov-proj}). More precisely, for
the first-order Lie--Trotter splitting we solve the equations
\begin{align}
\partial_{t}f & =-P_{\overline{V}}\left(v\cdot\nabla_{x}f-E(f)\cdot\nabla_{v}f\right),\label{eq:split-I}\\
\partial_{t}f & =P_{\overline{V}}P_{\overline{X}}\left(v\cdot\nabla_{x}f-E(f)\cdot\nabla_{v}f\right)\label{eq:split-II}\\
\partial_{t}f & =-P_{\overline{X}}\left(v\cdot\nabla_{x}f-E(f)\cdot\nabla_{v}f\right)\label{eq:split-III}
\end{align}
one after the other. In the following discussion we
consider just the first-order Lie--Trotter splitting algorithm with step size $\tau$.
The extension to second order (Strang splitting)
is nearly straightforward and is presented in Section \ref{subsec:Second-order}.

We assume that the initial value for the algorithm is given in the
 form
\[
f(0,x,v)=\sum_{i,j}X_{i}^{0}(x)S_{ij}^{0}V_{j}^{0}(v).
\]
First, let us consider equation (\ref{eq:split-I}). Since the set
$\{V_{j}\colon j=1,\dots,r\}$ forms an orthonormal basis of $\overline{V}$
(for each $t$), we have for the approximate particle density (again denoted by $f$)
\begin{equation}
f(t,x,v)=\sum_{j}K_{j}(t,x)V_{j}(t,v),\qquad\text{with}\quad K_{j}(t,x)=\sum_{i}X_{i}(t,x)S_{ij}(t),\label{eq:Kdef}
\end{equation}
where $K_{j}(t,x)$ is the coefficient of $V_{j}$ in the corresponding
basis expansion. We duly note that $K_{j}$ is a function of $x$,
but not of $v$. We then rewrite equation (\ref{eq:split-I}) as
follows
\begin{equation} \label{eq:pre-evol-K}
\begin{aligned}
&\sum_{j}\partial_{t}K_{j}(t,x)V_{j}(t,v)+\sum_{j}K_{j}(t,x)\partial_{t}V_{j}(t,v)
\\
&\hskip 2cm
=-\sum_{j}\bigl\la V_{j}(t,\cdot),v\mapsto v\cdot\nabla_{x}f(t,x,v)-E(f)(x)\cdot\nabla_{v}f(t,x,v)\bigr\ra_v V_{j}(t,v).
\end{aligned} \end{equation} 
    \change{Note that the electric field is self-consistently determined according to equation (\ref{eq:vlasov-poisson}). In practice this is done by solving the Poisson problem $-\Delta \phi = \rho(f)+1$ for the potential $\phi$. The electric field is then determined according to $E = -\nabla \phi$. We discuss how the (electron) charge density $\rho(f)=-\int f(t,x,v)\,\mathrm{d}v$ is computed in terms of the low-rank approximation after stating the corresponding evolution equations (for each step in the splitting algorithm).}

The solution of equation (\ref{eq:pre-evol-K}) is given by $V_{j}(t,v)=V_{j}(0,v)=V^{0}_j(v)$
and
\begin{equation}
\partial_t{K}_{j}(t,x)=-\sum_{l}c_{jl}^{1}\cdot\nabla_{x}K_{l}(t,x)+\sum_{l}c_{jl}^{2}\cdot E(K)(t,x)K_{l}(t,x)\label{eq:evolution-K}
\end{equation}
with
\[
c_{jl}^{1}=\int_{\Omega_v} vV_{j}^{0}V_{l}^{0}\,\mathrm{d}v,\qquad c_{jl}^{2}=\int_{\Omega_v} V_{j}^{0}(\nabla_{v}V_{l}^{0})\,\mathrm{d}v.
\]
The latter is obtained by equating coefficients in the basis expansion.
A very useful property of this splitting step is that we have to
update only the $K_{j}$, but not the $V_{j}$. We further note that
$c_{jl}^{1}$
and $c_{jl}^{2}$ are vectors in $\mathbb{R}^d$, since $v$ and $\nabla_{v}V_l^0$ are vector quantities. We denote the $m$th component of $c_{jl}^k$ by $c_{jl}^{k;x_m}$.

In the equation derived above we have written $E(K)$ to denote that the
electric field only depends on $K$ during the current step of the
algorithm. We will now consider this point in more detail. The electric
field is calculated from the electric charge density
\[
\rho(f)(t,x)=-\int_{\Omega_v} f(t,x,v)\,\mathrm{d}v
\]
\change{as described above}. Once the charge
$\rho$ is specified, the electric field is uniquely determined. For
the purpose of solving equation (\ref{eq:evolution-K}) we have the charge
\begin{align*}
\rho(t,x) & =-\sum_{j}K_{j}(t,x)\int_{\Omega_v} V_{j}^{0}(v)\,\mathrm{d}v\\
 & =-\sum_{j}K_{j}(t,x)\rho(V_{j}^{0}).
\end{align*}
Thus, the electric field only depends on $K=(K_{1},\dots,K_{r})$,
which explains our notation $E(K)$.

Both solving equation (\ref{eq:evolution-K}) and determining the
electric field are problems posed in a $d$-dimensional (as opposed
to $2d$-dimensional) space. Thus, we proceed by integrating equation
(\ref{eq:evolution-K}) with initial value
\[
K_{j}(0,x)=\sum_{i}X_{i}^{0}(x)S_{ij}^{0}
\]
until time $\tau$ to obtain $K_{j}^{1}(x)=K_{j}(\tau,x)$. However,
this is not sufficient as the $K_{j}^{1}$ are not necessarily orthogonal
(a requirement of our low-rank representation). Fortunately, this
is easily remedied by performing a QR decomposition
\[
K_{j}^{1}(x)=\sum_{i}X_{i}^{1}(x)\widehat S_{ij}^{1}
\]
to obtain orthonormal $X_{i}^{1}$ and the matrix $\widehat S_{ij}^{1}$. Once
a space discretization has been introduced, this QR decomposition
can be simply computed by using an appropriate function from a software
package such as LAPACK. However, from a mathematical point of view,
the continuous dependence on $x$ causes no issues. For example, the
modified Gram-Schmidt process works just as well in the continuous
formulation considered here.

Second, we proceed in a similar way for equation (\ref{eq:split-II}).
In this case both $V_{j}^{0}$ and $X_{i}^{1}$ are unchanged and
only $S_{ij}$ is updated. The corresponding evolution equation (which runs backward in time) is
given by
\begin{align}
\partial_t S_{ij}(t) & =\bigl\la X_{i}^{1}(x)V_{j}^{0}(v),(v\cdot\nabla_{x}-E(S)(t,x)\cdot\nabla_{v})\sum_{k,l} X_{k}^{1}(x)S_{kl}(t)V_{l}^{0}(v)\bigr\ra_{x,v}\nonumber \\
 & =\sum_{k,l}\left(c_{jl}^{1}\cdot d_{ik}^{2}-c_{jl}^{2}\cdot d_{ik}^{1}[E(S(t))]\right)S_{kl}(t)\label{eq:evolution-S}
\end{align}
with
\[
d_{ik}^{1}[E]=\int_{\Omega_x} X_{i}^{1}EX_{k}^{1}\,\mathrm{d}x,\qquad\quad d_{ik}^{2}=\int_{\Omega_x} X_{i}^{1}(\nabla_{x}X_{k}^{1})\,\mathrm{d}x,
\]
and $E(S(t))$ denotes the electric field corresponding to the charge density
\[
\rho(t,x)  = -\sum_{i,j}X_i^1(x) S_{ij}(t)\rho(V_{j}^{0}).
\]
Note that in this case the evolution equation depends neither on $x$
nor on $v$. Since $E$ and $\nabla_{x}X_k^1$ are vector quantities, we
have coefficient vectors $d_{ik}^{1}$ and $d_{ik}^{2}$ in $\mathbb{R}^d$. 
We now integrate equation (\ref{eq:evolution-S}) with initial value
$S_{ij}(0)=\widehat S_{ij}^{1}$ until time $\tau$ obtain $\widetilde S_{ij}^{0}=S_{ij}(\tau)$.
This completes the second step of the algorithm.

Finally, we consider equation (\ref{eq:split-III}). Similar to the
first step we have
\[
f(t,x,v)\approx\sum_{i}X_{i}(t,x)L_{i}(t,v),\qquad\text{with}\quad L_{i}(t,v)=\sum_{j}S_{ij}(t)V_{j}(t,v).
\]
As before, it is easy to show that this time the $X_{i}$ remain constant during
that step. Thus, the $L_{j}$ satisfy the following evolution equation
\begin{align}
\partial_t{L}_{i}(t,v) & =-\Bigl\la X_{j}^{1},(v\cdot\nabla_{x}-E(L)(x,t)\cdot\nabla_{v})\sum_{k}X_{k}^{1}L_{k}(t,v)\Bigr\ra_x \nonumber \\
 & =\sum_{k}d_{ik}^{1}[E(L(t,\cdot))]\cdot\nabla_{v}L_{k}(t,v)-\sum_{k}(d_{ik}^{2}\cdot v)L_{k}(t,v),\label{eq:evolution-L}
\end{align}
where $E(L(t,\cdot))$ denotes the electric field that corresponds to the charge density
$$
\rho(t,x)  = -\sum_{i}X_i^1(x) \rho(L_{i}(t,\cdot)).
$$
We then integrate equation (\ref{eq:evolution-L}) with initial value
\[
L_{i}(0,v)=\sum_{j}\widetilde S_{ij}^{0}V_{j}^{0}(v)
\]
up to time $\tau$ to obtain $L_{i}^{1}(v)=L_{i}(\tau,v)$. Since,
in general, the $L_{i}^{1}$ are not orthogonal we have to perform
a QR decomposition
\[
L_{i}^{1}(v)=\sum_{j}S_{ij}^{1}V_{j}^{1}(v)
\]
to obtain  orthonormal functions $V_{j}^{1}$ and the matrix $S_{ij}^{1}$. Finally, the output of our
Lie splitting algorithm is
\[
f(\tau,x,v)\approx\sum_{i,j}X_{i}^{1}(x)S_{ij}^{1}V_{j}^{1}(v).
\]

To render this algorithm into a numerical scheme that can be implemented
on a computer, we have to discretize in space. In principle, it is
then possible to integrate equations (\ref{eq:evolution-K}), (\ref{eq:evolution-S}),
and (\ref{eq:evolution-L}) using an arbitrary time stepping method.
However,  equations (\ref{eq:evolution-K}) and (\ref{eq:evolution-L})
share many common features with the original Vlasov system. Therefore, employing
a numerical method tailored to the present situation can result in
significant performance improvements. More precisely, we will consider
fast Fourier based techniques (Section \ref{subsec:FFT}) and a semi-Lagrangian
approach (Section \ref{subsec:Semi-Lagrangian-method}). Both of these
schemes remove the CFL condition from equations (\ref{eq:evolution-K})
and (\ref{eq:evolution-L}). Similar methods have been extensively
used for direct Eulerian simulation of the Vlasov equation (which
is the primary motivation to extend them to the low-rank algorithm
proposed in this work).

Before proceeding, however, let us briefly discuss the complexity
of our algorithm. For that purpose we choose to discretize both $x$
and $v$ using $n$ grid points in every coordinate direction. The storage cost of our numerical
method is then $\mathcal{O}(n^{d}r)$, compared to $\mathcal{O}(n^{2d})$
for the Eulerian approach. Now, let us assume that the cost of solving
the evolution equations is proportional to the degrees of freedom
(which is usually true for semi-Lagrangian schemes and true up to
a logarithm for spectral methods). In this case our low-rank algorithm
requires $\mathcal{O}(n^{d}r)$ arithmetic operations for the evolution equations and $O(n^dr^2)$ operations for computing the integrals, compared
to $\mathcal{O}(n^{2d})$ for a direct simulation.

\textbf{}

\subsection{Spectral method\label{subsec:FFT}}

In the framework of the classic splitting algorithm for the Vlasov\textendash Poisson
equation (which goes back to \cite{cheng1976}), spectral methods
have long been considered a viable approach. This is mainly due to
the fact that once the splitting has been performed, the remaining
advection equations can be solved efficiently by employing fast Fourier
transforms. Equations (\ref{eq:evolution-K}) and (\ref{eq:evolution-L})
are somewhat more complicated, compared to the advection equations
resulting from a direct splitting of the Vlasov\textendash Poisson
equation, in that a nonlinear term is added. In addition, for equation
(\ref{eq:evolution-L}) even the speed of the advection depends on
the electric field. However, since the advection coefficients are
scalar quantities it is still possible to apply Fourier techniques,
as we will show below.

We first consider equation (\ref{eq:evolution-K}). Performing a Fourier
transform (denoted by $\mathcal{F}$) in $x$ yields
\[
\partial_t{\hat{K}}_{j}(t,k)=-\sum_{\beta}\left(c_{jl}^{1}\cdot ik\right)\hat{K}_{l}(t,k)+\sum_{l}c_{jl}^{2}\cdot\mathcal{F}\left(E(K)(\cdot)K_{l}(t,\cdot)\right).
\]
Forming the vectors $\hat{K}=(\hat{K}_{1},\dots,\hat{K}_{r})$ and
$K=(K_{1},\dots,K_{r})$ and defining $F$ appropriately, we can write
this in matrix notation
\begin{equation}
\dot{\hat{K}}(t,k)=A(k)\hat{K}(t,k)+F(\hat{K}(t,\cdot))(k).\label{eq:A+F}
\end{equation}
The linear part is the only source of stiffness. For each $k$ we
have $A(k)\in\mathbb{R}^{r\times r}$ which is a small enough matrix
to be handled by direct methods. Thus, an exponential integrator is
easily able to integrate the linear part exactly (which removes the
CFL condition). For a first order scheme the exponential Euler method
is sufficient
\[
\hat{K}(\tau,k)\approx\mathrm{e}^{\tau A(k)}\hat{K}(0,k)+\tau\varphi_{1}(\tau A(k))F(\hat{K}(0,\cdot))(k),
\]
where $\varphi_{1}(z)=(\mathrm{e}^{z}-1)/z$ is an entire function.
It should be noted that higher order methods have been derived as
well (see, for example, the review article \cite{Hochbruck2010}).

Now, let us turn our attention to equation (\ref{eq:evolution-L}).
Before this equation can be made amenable to Fourier techniques,
a further approximation has to be performed. Specifically, we freeze
the electric field at the beginning of the time step. That is, instead
of equation (\ref{eq:evolution-L}) we solve the approximation
\[
\partial_t{L}_{i}(t,v)=\sum_{k}d_{ik}^{1}\cdot\nabla_{v}L_{k}(t,v)-\sum_{k}(d_{ik}^{2}\cdot v)L_{k}(t,v)
\]
with $d_{ik}^{1}=d_{ik}^{1}(E(L(0,\cdot)))$. This still yields a
first order approximation. We will show in Section \ref{subsec:Second-order}
how this technique can be improved to second order. Now, the advection
speed is independent of $L$ and we are able to apply a Fourier transform
in $v$ which yields
\[
\partial_t{\hat{L}}_{i}(t,k)=\sum_{k}\left(d_{ik}^{1}\cdot ik\right)\hat{L}_{k}(t,k)-\sum_{k}d_{ik}^{2}\cdot\mathcal{F}\left(vL_{k}(t,\cdot)\right).
\]
This is in the form of equation (\ref{eq:A+F}) and can thus be handled
by an exponential integrator in the same way as is explained above.

\subsection{Semi-Lagrangian method\label{subsec:Semi-Lagrangian-method}}

Semi-Lagrangian methods are widely employed to solve the Vlasov equation.
To a large part this is due to the fact that by performing a splitting
in the different spatial and velocity directions, only one-dimensional
advection equations have to be solved. The projection operation necessary
for these class of methods can then be implemented easily and efficiently.
In the present section our goal is to show that an efficient semi-Lagrangian
scheme can be derived for equations (\ref{eq:evolution-K}) and (\ref{eq:evolution-L}).

We start with equation (\ref{eq:evolution-K}) and use the vector
notation $K=(K_{1},\dots,K_{r})$. Then,
\begin{equation}
\partial_t{K}(t,x)=-c^{1;x_{1}}\partial_{x_{1}}K(t,x)-c^{1;x_{2}}\partial_{x_{2}}K(t,x)-c^{1;x_{3}}\partial_{x_{3}}K(t,x)+\left(c^{2}\cdot E(K)(t,x)\right)K(t,x).\label{eq:K-matrix}
\end{equation}
Applying a first order Lie splitting to equation (\ref{eq:K-matrix}),
with initial value $K(0,x)=K^{0}(x)$, gives 
\[
{K}(\tau,\cdot)\approx\mathrm{e}^{-\tau c^{1;x_{1}}\partial_{x_{1}}}\mathrm{e}^{-\tau c^{1;x_{2}}\partial_{x_{2}}}\mathrm{e}^{-\tau c^{1;x_{3}}\partial_{x_{3}}}(K^{0}+\tau(c^{2}\cdot E(K^{0}))K^{0}).
\]
This can be extended trivially to second or higher order. The crucial
part is the computation of
\[
M(t,x)=\mathrm{e}^{-tc^{1;x_{1}}\partial_{x_{1}}}M(0,x)
\]
which is equivalent to the partial differential equation
\[
\partial_t{M}(t,x)=-c^{1;x_{1}}\partial_{x_{1}}M(t,x).
\]
Since $c^{1;x_{1}}$ is symmetric, there exists an orthogonal matrix
$T$ such that $Tc^{1;x_{1}}T^{T}=D$, where $D$ is a diagonal matrix.
All these computations can be done efficiently as $c^{1;x_{1}}\in\mathbb{R}^{r\times r}$
(i.e.~we are dealing with small matrices). We now change variables
to $\overline{M}=TM$ and obtain
\[
\partial_t{\overline{M}}_{j}(t,x)=-D_{jj}\partial_{x_{1}}\overline{M}(t,x).
\]
This is a one-dimensional advection equation with constants coefficients
and can thus be treated with an arbitrary semi-Lagrangian method.
Once $\overline{M}(\tau,x)$ is obtained the coordinate transformation
is undone such that we get $M_{j}(\tau,x)=T^{T}\overline{M}_{j}(\tau,x)$.
Exactly the same procedure is then applied to the advection in the
$x_{2}$ and $x_{3}$ direction.

For equation (\ref{eq:evolution-L}) we proceed in a similar way.
The main difference here is that we first freeze $d_{ik}^{1}=d_{ik}^{1}(E(L(0,\cdot))$
at the beginning of the step (as is explained in more detail in the
previous section). This leaves us with 
\[
\partial_t{L}_{i}(t,v)=\sum_{k}d_{ik}^{1}\cdot\nabla_{v}L_{k}(t,v)-\sum_{k}(d_{ik}^{2}\cdot v)L_{k}(t,v).
\]
Since this equation is (almost) identical to equation (\ref{eq:K-matrix})
and $d^{1;x_{1}}$, $d^{1;x_{2}}$, and $d^{1;x_{3}}$ are symmetric,
we can apply the same algorithm.

\subsection{Second-order integrator\label{subsec:Second-order}}

So far we have only considered the first order Lie splitting. However,
in principle it is simple to employ the second order Strang splitting.
The main obstacle in this case is that for both the Fourier approach
(introduced in section \ref{subsec:FFT}) and the semi-Lagrangian
approach (introduced in section \ref{subsec:Semi-Lagrangian-method})
it is essential that a further approximation is made when solving
equation (\ref{eq:evolution-L}). More specifically, the electric
field $E$ is assumed fixed during that step. Thus, a naive application
of Strang splitting would still only result in a first order method.
However, in \cite{einkemmer1306,einkemmer1309} a technique has been
developed that can overcome this difficulty for methods of arbitrary
order.

A detailed account of the second order accurate scheme is given in
Algorithm \ref{alg:strang-splitting}. From there it should also be
clear that no additional ingredients, compared to the first order
scheme described in section \ref{sec:x-v decomposition}, are needed.
The main point is that we use an embedded Lie scheme of step size
$\tau/2$ in order to determine an electric field $E^{1/2}$ such
that $E^{1/2}(x)\approx E(f(\tfrac{\tau}{2},\cdot))(x)$. Using $E^{1/2}$
in the middle step of the Strang splitting, which approximates equation
(\ref{eq:evolution-L}), yields a numerical method that is almost
symmetric. Then composing two Lie steps in opposite order yields a
second order method. 

\begin{algorithm}
\begin{algorithmic}[1]
\Require{$X_i^0$, $S_{ij}^0$, $V_j^0$ (such that $f(0,x,v)\approx \sum_{i,j} X_i^0(x)S_{ij}^0 V_j^0(v))$}
\Ensure{$X_i^2$, $S_{ij}^5$, $V_j^1$ (such that $f(\tau,x,v)\approx \sum_{i,j} X_i^2(x)S_{ij}^5 V_j^1(v))$}
\State{Solve equation (\ref{eq:evolution-K}) with initial value $\sum_i X_i^0 S_{ij}^0$ up to time $\tau/2$ to obtain $K_j^1$.}
\State{Perform a QR decomposition of $K_j^1$ to obtain $X_i^1$ and $S_{ij}^1$.}
\State{Solve equation (\ref{eq:evolution-S}) with initial value $S_{ij}^1$ up to time $\tau/2$ to obtain $S_{ij}^2$.}
\State{Compute the electric field $E$ using $X_i^1$, $S_{ij}^2$, and $V_j^0$.}
\State{Compute $d^1_{ij}$ using $E$.}
\State{Solve equation (\ref{eq:evolution-L}) with fixed $d^1_{ij}$ and initial value $\sum_j S_{ij}^2 V_j^0$ up to time $\tau/2$ to obtain $L_{i}^{1/2}$.}
\State{Compute the electric field $E^{1/2}$ using $L_i^{1/2}$ and $X_i^1$.}
\State{Compute $d^1_{ij}$ using $E^{1/2}$.}
\State{Solve equation (\ref{eq:evolution-L}) with the fixed $d^1_{ij}$ and initial value $\sum_j S_{ij}^2 V_j^0$ up to time $\tau$ to obtain $L_i^1$.}
\State{Perform a QR decomposition of $L_i^1$ to obtain $V_j^1$ and $S_{ij}^3$.}
\State{Solve equation (\ref{eq:evolution-S}) with initial value $S_{ij}^3$ up to time $\tau/2$ to obtain $S_{ij}^4$.}
\State{Solve equation (\ref{eq:evolution-K}) with initial value $\sum_i X_i^1 S_{ij}^4$ up to time $\tau/2$ to obtain $K_j^2$.}
\State{Perform a QR decomposition of $K_j^2$ to obtain $X_i^2$ and $S_{ij}^5$.}
\end{algorithmic}\caption{A second-order accurate low-rank algorithm for the Vlasov\textendash Poisson
equation. \label{alg:strang-splitting}}
\end{algorithm}

\section{A hierarchical low-rank algorithm\label{sec:The-hierarchical-algorithm}}

The algorithm proposed in the previous section reduces a $2d$-dimensional
problem to a number of $d$-dimensional problems. However, if the
solution under consideration also admits a low-rank structure in $X_{i}(t,x)$
and $V_{j}(t,v)$, it is possible to further reduce the dimensionality
of the equations that need to be solved. To derive such an algorithm
is the content of the present section. We will first
consider the four-dimensional case ($d=2$). This is sufficient to elucidate
the important aspects of the algorithm. The algorithm for the six-dimensional
($d=3$) case is stated in Section \ref{subsec:h6d}.

\subsection{The hierarchical low-rank algorithm in four dimensions\label{subsec:h4d}}

Once again we start from the low-rank approximation
\begin{equation}
f(t,x,v)\approx \sum_{i,j}X_{i}(t,x)S_{ij}(t)V_{j}(t,v),\label{eq:f-rep-hierarch}
\end{equation}
where $x=(x_{1},x_{2})\in\Omega_{x}\subset\mathbb{R}^{2}$, $v=(v_{1},v_{2})\in\Omega_{v}\subset\mathbb{R}^{2}$,
and $S\in\mathbb{R}^{r\times r}$. In addition, we will now restrict $X_i$ and $V_j$ to
the low-rank representations
\begin{equation}
X_{i}(t,x)=\sum_{\alpha,\beta}X_{1\alpha}(t,x_{1})C_{i\alpha\beta}(t)X_{2\beta}(t,x_{2})\label{eq:hX-lowrank}
\end{equation}
and
\begin{equation}
V_{j}(t,v)=\sum_{\alpha,\beta}V_{1\alpha}(t,v_{1})D_{j\alpha\beta}(t)V_{2\beta}(t,v_{2}),\label{eq:hV-lowrank}
\end{equation}
where $C\in\mathbb{R}^{r\times r_{x}\times r_{x}}$, $D\in\mathbb{R}^{r\times r_{v}\times r_{v}}$
and $r_{x}$ and $r_{v}$ is the rank in the $x$- and $v$-direction,
respectively. It is, however, entirely reasonable for the Vlasov equation
that, for example, a further low-rank structure is present in velocity
space but the same is not true in physical space. \change{As an example, consider
a plasma system, where only weak kinetic effects are present. Traditionally, fluid models (such
as magnetohydrodynamics and its extensions) have been used to model such systems. However, this
approach neglects kinetic effects altogether. As an alternative the low-rank approach in this
paper could be employed. Since the system is still close to thermodynamic equilibrium, a further
low-rank approximation can be employed in $v$. On the other hand, the dynamics in $x$ might not
be accessible to such an approximation.}
In this case only the hierarchical low-rank representation given in (\ref{eq:hV-lowrank})
would be used. This would leave equation (\ref{eq:evolution-K}) precisely
as stated in Section \ref{sec:x-v decomposition}, while equation
(\ref{eq:evolution-L}) is replaced by the algorithm developed in
this section. 

Fortunately, equations (\ref{eq:evolution-K}) and (\ref{eq:evolution-L})
are sufficiently similar that we only have to derive (and perhaps
more importantly, implement) the algorithm once. Thus, we will consider
the following equation
\[
\partial_t{K}_{j}(t,x)=-\sum_{l}c_{jl}^{1}\cdot\nabla_{x}K_{l}(t,x)-\sum_{l}c_{jl}^{2}\cdot F(K)(x)K_{l}(t,x).
\]
For solving equation (\ref{eq:evolution-K}) we simply set $F(K)=-E(K)$.
If we replace $c^{1}$ and $c^{2}$ by $d^{1}$ and $d^{2}$, respectively, replace the variable $x$ by $v$\change{, reverse the sign of the equation}, and set
$F(K)(v)=v$, then we obtain precisely equation (\ref{eq:evolution-L}).

In the following we transfer the algorithm proposed for time integration of Tucker
tensors in \cite{Lubich2017} to the current setting. The reader is encouraged to first look up \cite{Lubich2017} before entering into the thicket of formulas presented in the following.
Compared to the generic case in \cite{Lubich2017}, there are
important simplifications for the problem under consideration, which
even influence the representation of our low-rank approximation.
In fact, a Tucker representation for the order-three tensor $X_{j}(t,x_{1},x_{2})$ (for fixed $t$)
would be given by
\[
X_{j}(t,x)=\sum_{l,\alpha,\beta}R_{l\alpha\beta}(t)A_{lj}(t)X_{1\alpha}(t,x_{1})X_{2\beta}(t,x_{2})
\]
and not by equation (\ref{eq:hX-lowrank}). The algorithm would then
proceed by updating $A$, $X_{1}$, $X_{2}$, and finally the 
core tensor $R$. However, in the present setting $A$ is artificial
in the sense that we can always go back to the representation in equation
(\ref{eq:hX-lowrank}) by setting
\[
C_{j\alpha\beta}(t)=\sum_{l}R_{l\alpha\beta}(t)A_{lj}(t).
\]
If we want to express this in the Tucker format, we can insert an
identity matrix
\[
X_{j}(t,x)=\sum_{l,\alpha,\beta}C_{l\alpha\beta}(t)\delta_{lj}X_{\alpha}^{1}(x_{1})X_{\beta}^{2}(x_{2}).
\]
It is straightforward to show that if we apply the projection splitting
algorithm in \cite{Lubich2017} to a Tucker tensor in this form, the
first step leaves both $C$ and $A$ (in this case the identity) unchanged.
Thus, we can proceed with the algorithm and work directly in the representation
given by (\ref{eq:hX-lowrank}). In the following we will consider
the first-order Lie--Trotter splitting. We assume that the initial value is
given in the form
\begin{equation} \label{eq:ivp-4d-f}
f(0,x,v)=\sum_{i,j}X_{i}^{0}(x)S_{ij}^{0}V_{j}^{0}(v)
\end{equation}
with
\begin{equation}  \label{eq:ivp-4d-X}
X_{i}^{0}(x)=\sum_{\alpha,\beta}X_{1\alpha}^{0}(x_{1})G_{i\alpha\beta}^{0}X_{2\beta}^{0}(x_{2}).
\end{equation}
Since the evolution equation is formulated for $K_j=\sum_i X_i S_{ij} $, we then compute
\[
K_{j}^{0}(x)=\sum_{\alpha,\beta}X_{1\alpha}^{0}(x_{1})C_{j\alpha\beta}^{0}X_{2\beta}^{0}(x_{2}),
\]
where
\[
C_{j\alpha\beta}^{0}=\sum_{i}G_{i\alpha\beta}^{0}S_{ij}^{0}.
\]

\textbf{Step 1: } The first step of the algorithm
updates $X_{1}$ and $C$. First, we perform a QR decomposition
\[
C_{j\alpha\beta}^{0}=\sum_{\xi}Q_{j\xi\beta}R_{\alpha\xi}^{0}
\]
and set
\begin{equation}
W_{j\xi}(x_{2})=\sum_{\beta}Q_{j\xi\beta}X_{2\beta}^{0}(x_{2}).\label{eq:W-step1}
\end{equation}
In an actual implementation, the QR decomposition can be computed
by defining $A\in\mathbb{R}^{(r\cdot r_{x})\times r_{x}}$ such that
$A_{(j+\beta r)\alpha}=C_{j\alpha\beta}^{0}$ and then computing
the QR decomposition of $A$. Now, let us define
\[
M_{\alpha}(t,x_{1})=\sum_{\xi}X_{1\xi}(t,x_{1})R_{\xi\alpha}(t),
\]
for which we give an evolution equation. These $M_{\alpha}$ in fact
play the same role as the $K_{j}$ did in the algorithm of Section
\ref{sec:x-v decomposition}. The equation derived is as follows:
\begin{align}
\partial_t{M}_{\alpha}(t,x_{1}) & =-\sum_{i}\Bigl\la W_{i\alpha},\sum_{k}\bigl(c_{ik}^{1}\cdot\nabla_{x}+c_{ik}^{2}\cdot F(M)(x_{1},\cdot)\bigr)\sum_{\xi}M_{\xi}(t,x_{1})W_{k\xi}\Bigr\ra_{x_2}\nonumber \\
 & =-\sum_{\xi}a_{\alpha\xi}^{1}\partial_{x_{1}}M_{\xi}(t,x_{1})-\sum_{\xi}\bigl(a_{\alpha\xi}^{2}+a_{\alpha\xi}^{3}(x_{1})\bigr) M_{\xi}(t,x_{1}),\label{eq:evol-M}
\end{align}
where 
\[
a_{\alpha\xi}^{1}=\sum_{i,k}c_{ik}^{1;x_{1}}\la W_{i\alpha},W_{k\xi}\ra_{x_2},\qquad a_{\alpha\xi}^{2}=\sum_{i,k}c_{ik}^{1;x_{2}}\la W_{i\alpha},\partial_{x_{2}}W_{k\xi}\ra_{x_2},
\]
and
\[
a_{\alpha\xi}^{3}(x_{1})=\sum_{i,k}c_{ik}^{2}\cdot\la W_{i\alpha},F(M)(x_{1},\cdot)W_{k\xi}\ra_{x_2}.
\]
Note that $a^{1}$, $a^{2}$, and $a^{3}$ can be easily expressed
in terms of $Q$ and inner products over functions involving $X_{2\beta}$.
For $a_{\alpha\beta}^{1}$ the orthonormality relation of the $X_{2\beta}$
can be used to obtain a simpler formula that can be computed without
performing any integrals:
\[
a_{\alpha\xi}^{1}=\sum_{i,k,\beta}c_{ik}^{1;x}Q_{i\alpha\beta}Q_{k\xi\beta}.
\]
For $a^{2}$ we have
\change{\[
    a_{\alpha\xi}^{2}=\sum_{k,\beta}\Bigl( \sum_i c_{ik}^{1;x_{2}}Q_{i\alpha\beta}\Bigr)\Bigl(\sum_{\eta} Q_{k\xi\eta}\la X_{2\beta}^{0},\partial_{x_{2}}X_{2\eta}^{0}\ra_{x_2}\Bigr)
\]}
which requires $\mathcal{O}(n\mathcal{R}^{2})$, with $\mathcal{R}=\text{max}\left(r,r_{x},r_{v}\right)$,
arithmetic operations to compute the integrals and then \change{$\mathcal{O}(\mathcal{R}^{4})$}
operations to obtain the entries $a_{\alpha\xi}^{2}$. A similar
formula can be given for $a^{3}$ but we postpone this discussion
until later (where we will also discuss the low-rank representation
of the electric field in more detail).

Now, we solve equation (\ref{eq:evol-M}) with initial value
\[
M_{\alpha}(0,x_{1})=\sum_{\xi}X_{1\xi}^{0}(x_{1})R_{\xi\alpha}^{0}
\]
until time $\tau$ and obtain $M_{\alpha}^{1}(x_{1})=M_{\alpha}(\tau,x_{1})$.
Then a QR factorization
\[
M_{\alpha}^{1}(x_{1})=\sum_{\xi}X_{1\xi}^{1}(x_{1})R_{\xi\alpha}^{1}
\]
is performed to obtain orthonormal functions $X_{1\xi}^{1}$ and the matrix $R^{1}$.

Now, an evolution equation for $R$ is derived to become
\begin{align}
\dot{R}_{\xi\alpha}(t) & =\sum_{i}\Bigl\la X_{1\xi}^{1}W_{i\alpha},\sum_{k}\bigl(c_{ik}^{1}\cdot\nabla_{x}+c_{ik}^{2}\cdot F(R)\bigr)\sum_{\xi',\alpha'}X_{1\xi'}^{1}R_{\xi'\alpha'}(t)W_{k\alpha'}\Bigr\ra_x\nonumber \\
 & =\sum_{\xi',\alpha'}b_{\xi\xi'}^{1}R_{\xi'\alpha'}(t)a_{\alpha\alpha'}^{1}+\sum_{\alpha'}R_{\xi\alpha'}(t)a_{\alpha\alpha'}^{2}+\sum_{\xi',\alpha'}B_{\xi\alpha\xi'\alpha'}R_{\xi'\alpha'}(t)\label{eq:evol-R-step1}
\end{align}
with
\[
b_{\xi\xi'}^{1}=\la X_{1\xi}^{1},\partial_{x_{1}}X_{1\xi'}^{1}\ra_{x_1},\qquad B_{\xi\alpha\xi'\alpha'}=\sum_{i,k}c_{ik}^{2}\cdot\la X_{1\xi}^{1}W_{i\alpha},F(R)X_{1\xi'}^{1}W_{k\alpha'}\ra_{x}.
\]
Let us note that $B_{\xi\alpha\xi'\alpha'}$ is defined by a two-dimensional
integral. However, it is, of course, possible to also obtain a low-rank
approximation for $F(R)$. Only two cases are relevant here. First,
$F(x)=x$ which trivially is a low-rank approximation of rank $1$.
Second, $F(R)=-E(R)$. Then, the electric field $E$ is approximated
by
\[
E(x_{1},x_{2})=\sum_{\mu}E_{1\mu}(x_{1}) \circ E_{2\mu}(x_{2}),
\]
where the component-wise product is written as $\circ$ and 
$\mu=1,\dots,r_{E}$; here
$r_{E}$ is the rank used to represent the electric field. Then
\begin{equation}
B_{\xi\alpha\xi'\alpha'}=-\sum_{i,k,\mu}c_{ik}^{2} \,\la X_{1\xi}^{1},E_{1\mu}X_{1\xi'}^{1}\ra_{x_1}\circ \la W_{i\alpha},E_{2\mu}W_{k\alpha'}\ra_{x_2},\label{eq:E-lowrnak}
\end{equation}
which only requires the evaluation of one-dimensional integrals. Nevertheless,
evaluating equation (\ref{eq:E-lowrnak}) naively requires $\mathcal{O}(\mathcal{R}^{7})$
arithmetic operations (with $\mathcal{R}=\text{max}\left(r,r_{x},r_{E}\right)$),
in addition to evaluating the integrals. However, we can write
\[
B_{\xi\alpha\xi'\alpha'}=-\sum_{\mu}\la X_{1\xi}^{1},E_{1\mu}X_{1\xi'}^{1}\ra_{x_1}\cdot\sum_{i,k}c_{ik}^{2} \la W_{i\alpha},E_{2\mu}W_{k\alpha'}\ra_{x_2}
\]
and evaluate the sums from right to left, which only requires $\mathcal{O}(\mathcal{R}^{5})$
arithmetic operations.

Now, we integrate equation (\ref{eq:evol-R-step1}) with initial value
$R_{\xi\alpha}(0)=R_{\xi\alpha}^{1}$ until time $\tau$ to obtain
$R_{\xi\alpha}^{2}=R_{\xi\alpha}(\tau)$. From this we update
$C$ as 
\[
C_{j\alpha\beta}^{1}=\sum_{\xi}Q_{j\xi\beta}R_{\alpha\xi}^{2},
\]
which completes the first step of the algorithm. 

\textbf{Step 2: }
The second step proceeds in a similar manner, but updates $X_{2}$
and $C$. First, we perform a QR decomposition
\[
C_{i\alpha\beta}^{1}=\sum_{\eta}Q_{i\alpha\eta}R_{\beta\eta}^{0}.
\]
The tensor $Q$ and the matrix $R^0$ obtained here are different from those in the first step, but for ease of notation we do not use different symbols or extra superscripts.
We define
\[
N_{\beta}(t,x_{2})=\sum_{\eta}X_{2\eta}(t,x_{2})R_{\eta\beta}(t)
\]
for which we determine an evolution equation. These $N_{\beta}$  play the same role as the $L_{j}$ did in the algorithm of Section
\ref{sec:x-v decomposition}. The equation derived is as follows:
\begin{align}
\partial_t{N}_{\beta}(t,x_{2}) & =-\sum_{i,\alpha}Q_{i\alpha\beta}\Bigl\la X_{1\alpha}^{1},\sum_{k}\left(c_{ik}^{1}\cdot\nabla_{x}+c_{ik}^{2}\cdot F(N)\right)\sum_{\xi,\eta}X_{1\xi}^{1}Q_{k\xi\eta}N_{\eta}(t,x_{2})\Bigr\ra_{x_1} \nonumber \\
 & =-\sum_{\eta}h_{\beta\eta}^{1}\partial_{x_{2}}N_{\eta}(t,x_{2})-\sum_{\eta}\left(h_{\beta\eta}^{2}+h_{\beta\eta}^{3}(x_{2})\right)N_{\eta}(t,x_{2}),\label{eq:evol-N-step2}
\end{align}
where 
\[
h_{\beta\eta}^{1}=\sum_{i,k,\alpha}c_{ik}^{1;x_{2}}Q_{i\alpha\beta}Q_{k\alpha\eta},\qquad h_{\beta\eta}^{2}=\sum_{i,k,\alpha,\xi}Q_{i\alpha\beta}Q_{k\xi\eta}c_{ik}^{1;x_{1}}b_{\alpha\xi}^{1}
\]
and 
\[
h_{\beta\eta}^{3}(x_{2})=\sum_{i,k,\alpha,\xi}Q_{i\alpha\beta}Q_{k\xi\eta}c_{ik}^{2}\cdot\bigl\la X_{1\alpha}^{1},F(N)(\cdot,x_{2})X_{1\xi}^{1}\bigr\ra_{x_1}.
\]
Equation (\ref{eq:evol-N-step2}) is integrated with initial value
\[
N_{\beta}(0,x_{2})=\sum_{\eta}X_{2\eta}^{0}(x_{2})R_{\eta\beta}^{0}
\]
until time $\tau$. We then set $N_{\beta}^{1}(x_{2})=N_{\beta}(\tau,x_{2})$
and perform a QR decomposition
\[
N_{\beta}^{1}=\sum_{\eta}X_{2\eta}^{1}R_{\eta\beta}^{1}
\]
to obtain orthonormal functions $X_{2\eta}^{1}$ and the matrix $R^{1}$. Next, an evolution equation for
$R$ is obtained as follows:
\begin{align}
\dot{R}_{\eta\beta}(t) & =\sum_{i,\alpha}Q_{i\alpha\beta}\Bigl\la X_{1\alpha}^{1}X_{2\eta}^{1},\sum_{k}\left(c_{ik}^{1}\cdot\nabla_{x}+c_{ik}^{2}\cdot F(R)\right)\sum_{\alpha',\eta',\beta'}X_{1\alpha'}^{1}Q_{k\alpha'\beta'}R_{\eta'\beta'}X_{2\eta'}^{1}\Bigr\ra_{x_1,x_2}\nonumber \\
 & =\sum_{\eta',\beta'}d_{\eta\eta'}^{1}R_{\eta'\beta'}(t)h_{\beta\beta'}^{1}+\sum_{\beta'}R_{\eta\beta'}(t)h_{\beta\beta'}^{2}+\sum_{\eta',\beta'}G_{\eta\beta\eta'\beta'}R_{\eta'\beta'}(t)\label{eq:evol-R-step2}
\end{align}
with
\[
d_{\eta\eta'}^{1}=\la X_{2\eta}^{1},\partial_{x_{2}}X_{2\eta'}^{1}\ra_{x_2},\qquad G_{\eta\beta\eta'\beta'}=\sum_{i,k,\alpha,\alpha'}Q_{i\alpha\beta}Q_{k\alpha'\beta'}c_{ik}^{2}\cdot\la X_{1\alpha}^{1}X_{2\eta}^{1},F(R)X_{1\alpha'}^{1}X_{2\eta'}^{1}\ra_{x_1,x_2}.
\]
We then integrate equation (\ref{eq:evol-R-step2}) with initial value
$R^{1}$ up to time $\tau$ and set $R_{\eta\beta}^{2}=R_{\eta\beta}(\tau)$.
We then update the core tensor $C$ as 
\[
C_{i\alpha\beta}^{2}=\sum_{\eta}Q_{i\alpha\eta}R_{\beta\eta}^{2}.
\]
This completes the second step of the algorithm.

\textbf{Step 3: }In the last step we  update $C$ and the matrix $S$.  The tensor $C$ is first updated by the following evolution equation:
\begin{align}
\dot{C}_{i\alpha\beta}(t) & =-\Bigl\la X_{1\alpha}^{1}X_{2\beta}^{1},\sum_{k}\left(c_{ik}^{1}\cdot\nabla_{x}+c_{ik}^{2}\cdot F(C)\right)\sum_{\xi,\eta}X_{1\xi}^{1}C_{k\xi\eta}(t)X_{2\eta}^{1}\Bigr\ra_x\nonumber \\
 & =-\sum_{k,\xi}c_{ik}^{1;x_{1}}b_{\alpha\xi}^{1}C_{k\xi\beta}(t)-\sum_{k,\eta}c_{ik}^{1;x_{2}}d_{\beta\eta}^{1}C_{k\alpha\eta}(t)-\sum_{k,\xi,\eta}(c_{ik}^{2}\cdot e_{\alpha\beta\xi\eta})C_{k\xi\eta}(t),\label{eq:evol-C-step3}
\end{align}
where
\[
e_{\alpha\beta\xi\eta}=\la X_{1\alpha}X_{2\beta},F(C)X_{1\xi}X_{2\eta}\ra_x.
\]
We integrate equation (\ref{eq:evol-C-step3}) with initial value
$C^{2}$ until time $\tau$ and set $C_{i\alpha\beta}^{3}=C_{i\alpha\beta}(\tau)$.
Now, we have obtained the following representation
\[
K_{j}(\tau,x)\approx K_{j}^{1}(x)=\sum_{\alpha,\beta}X_{1\alpha}^{1}(x_{1})C_{j\alpha\beta}^{3}X_{2\beta}^{1}(x_{2}).
\]
This, however, is not yet sufficient, since in order to obtain an approximation
in the form given by equation (\ref{eq:f-rep-hierarch}), we have
to perform a QR decomposition. If this is done naively, the complexity
would scale as $n^{d}$, with $n$ the number of grid points per direction,
which is precisely what we want to avoid in the present setting. However,
since the $X_{1}^{1}$ and $X_{2}^{1}$ are already orthogonalized,
we can compute a QR decomposition without modifying $X_{1}^{1}$ or
$X_{2}^{1}$ and without evaluating any integrals. This is accomplished,
for example, by carrying out the modified Gram-Schmidt process. To
perform the Gram-Schmidt process, inner products and linear combinations
of the different $K_{j}^{1}$ are required. However, since by orthogonality
\begin{align}
\la K_{j}^{1},K_{k}^{1}\ra_x & =\sum_{\alpha\beta\xi\eta}C_{j\alpha\beta}^{3}C_{k\xi\eta}^{3}\la X_{1\alpha}^{1},X_{1\xi}^{1}\ra_{x_1}\la X_{2\beta}^{1},X_{2\eta}^{1}\ra_{x_2}\nonumber \\
 & =\sum_{\alpha,\beta}C_{j\alpha\beta}^{3}C_{k\alpha\beta}^{3}\label{eq:gs-inner}
\end{align}
and since
\begin{equation}
K_{j}^{1}+cK_{k}^{1}=\sum_{\alpha,\beta}X_{1\alpha}^{1}(x_{1})C_{\alpha\beta}^{+}X_{2\beta}^{1}(x_{2})\quad\text{ with }\quad  C_{\alpha\beta}^{+}=C_{j\alpha\beta}^{3}+cC_{k\alpha\beta}^{3},
\label{eq:gs-linear}
\end{equation}
this can be done efficiently, as only operations on $C$ have to be
performed. The corresponding procedure is shown in Algorithm \ref{alg:gram-schmidt} \change{and requires $\mathcal{R}^4$ arithmetic operations. It is applied to the input $C_{j\alpha\beta}^3$ and results in the output $C_{j\alpha\beta}^4$ and $S_{ij}^1$ (in an actual implementation this can be done in place, as is demonstrated by \ref{alg:gram-schmidt}). Thus, we finally obtain the low-rank representation
\[
f(0,x,v)=\sum_{i,j}X_{i}^{1}(x)S_{ij}^{1}V_{j}^{1}(v)
\]
with
\[
X_{i}^{1}(x)=\sum_{\alpha,\beta}X_{1\alpha}^{1}(x_{1})C_{i\alpha\beta}^{4}X_{2\beta}^{1}(x_{2}).
\]
This is precisely the same form in which the initial value is provided, see equations (\ref{eq:ivp-4d-f}) and (\ref{eq:ivp-4d-X}). We then can proceed to apply the corresponding procedure to update the quantities depending on velocity space.
}

\begin{algorithm}[h]
\begin{algorithmic}[1]
\Require{$C_{j\alpha\beta}$}
\Ensure{$S_{ij}$, $C_{j\alpha\beta}$}
\State{$S=0$}
\For{$j=1,\dots,r$}
  \For{$k=1,\dots,j-1$}
    \State{$S_{kj} = \sum_{\alpha \beta} C_{k\alpha \beta} C_{j \alpha \beta}$}
       \For{$\alpha$, $\beta$}
         \State{$C_{j\alpha\beta} = C_{j\alpha\beta} - S_{kj} C_{k\alpha\beta}$}
      \EndFor
    \EndFor
	\State{$S_{jj} = \sqrt{\sum_{\alpha,\beta} C_{j\alpha\beta}^2}$}
    \For{$\alpha$, $\beta$}
        \State{$C_{j\alpha\beta} = C_{j\alpha\beta}/S_{jj}$}
    \EndFor
\EndFor
\end{algorithmic}

\caption{The goal of the algorithm is to obtain $X_{i}$ and $S_{ij}$ such
that $\la X_{i},X_{k}\ra_x=\delta_{ik}$ and $K_{j}=\sum_{i}X_{i}S_{ij}$
with $X_{i}=\sum_{\alpha,\beta}X_{1\alpha}^{1}C_{i\alpha\beta}X_{2\beta}^{1}$.
The input is such that $K_{j}=\sum_{\alpha,\beta}X_{1\alpha}^{1}C_{j\alpha\beta}X_{2\beta}^{1}$
and consequently only $C$ is modified by the algorithm. This is accomplished
using the modified Gram-Schmidt procedure and equations (\ref{eq:gs-inner})
and (\ref{eq:gs-linear}). \label{alg:gram-schmidt}}
\end{algorithm}

\change{We have not yet discussed how to solve the evolution equations derived in this section. Equations (\ref{eq:evol-M}) and (\ref{eq:evol-N-step2}) (which correspond to an advection in $x_1$ and $x_2$, respectively) can be treated using the approach outlined in sections \ref{subsec:FFT} and \ref{subsec:Semi-Lagrangian-method}. That is, we can employ a spectral method or a semi-Lagrangian scheme. In fact, the present situation is simpler since the evolution equations only depend on a single variable. In fact, the further splitting into the different coordinate axis, as described in section \ref{subsec:Semi-Lagrangian-method} for the semi-Lagrangian approach, is not required in the present setting.
}

\change{Before proceeding let us discuss the complexity of the proposed hierarchical low-rank splitting scheme. Solving the evolution equations (\ref{eq:evol-M}) and (\ref{eq:evol-N-step2}) requires $\mathcal{O}(n\mathcal{R}^2)$ arithmetic operations and $\mathcal{O}(n\mathcal{R})$ storage, where $n$ is the number of grid points (per direction) and $\mathcal{R}=\max(r,r_x,r_v,r_E)$. In addition, we have to compute the coefficients which requires $\mathcal{O}(n\mathcal{R}^2+\mathcal{R}^4)$ arithmetic operations. The storage cost of the coefficients is $\mathcal{O}(n\mathcal{R}^2+\mathcal{R}^4)$ and is dominated by $a^3$ and $h^3$.}

\subsection{The hierarchical low-rank algorithm in six dimensions\label{subsec:h6d}}

In this section we extend the hierarchical low-rank
approximation to the six-dimensional ($d=3$) case. Since the derivation
is very similar we will only state the relevant evolution equations
here. We start with the low-rank representation
\begin{equation}
f(t,x,v)\approx \sum_{i,j}X_{i}(t,x)S_{ij}(t)V_{j}(t,v),\label{eq:f-rep-hierarch-1}
\end{equation}
where $x=(x_{1},x_{2},x_{3})\in\Omega_{x}\subset\mathbb{R}^{3}$,
$v=(v_{1},v_{2},v_{3})\in\Omega_{v}\subset\mathbb{R}^{3}$, and $S\in\mathbb{R}^{r\times r}$.
We now restrict to the low-rank representation 
\begin{equation}
X_{i}(t,x)=\sum_{\alpha,\beta,\gamma}C_{i\alpha\beta\gamma}(t)X_{1\alpha}(t,x_{1})X_{2\beta}(t,x_{2})X_{3\gamma}(t,x_{3})\label{eq:hX-lowrank-1}
\end{equation}
and
\begin{equation}
V_{j}(t,v)=\sum_{\alpha,\beta,\gamma}D_{j\alpha\beta\gamma}(t)V_{1\alpha}(t,v_{1})V_{2\beta}(t,v_{2})V_{3\gamma}(t,v_{3}),\label{eq:hV-lowrank-1}
\end{equation}
where $C\in\mathbb{R}^{r\times r_{x}\times r_{x}\times r_{x}}$, $D\in\mathbb{R}^{r\times r_{v}\times r_{v}\times r_{v}}$
and $r_{x}$ and $r_{v}$ is the rank in the $x$- and $v$-direction,
respectively. We then consider the evolution equation
\[
\partial_t{K}_{j}(t,x)=-\sum_{l}c_{jl}^{1}\cdot\nabla_{x}K_{l}(t,x)-\sum_{l}c_{jl}^{2}\cdot F(K)(x)K_{l}(t,x),
\]
which, depending on the choice of $F$, models either the update for
$X_{i}$ or $V_{j}$. Like in the previous section, we start
with the initial value
\[
K_{j}^{0}(x)=\sum_{\alpha,\beta,\gamma}C_{j\alpha\beta\gamma}^{0}X_{1\alpha}^{0}(x_{1})X_{2\beta}^{0}(x_{2})X_{3\gamma}^{0}(x_{3}).
\]
In the following we will divide the algorithm into four parts which
correspond to the update of $X_{1}$ and $C$, the update of $X_{2}$
and $C$, the update of $X_{3}$ and $C$, and finally an update of
$C$.

\textbf{Step 1: }We perform a QR decomposition
\[
C_{j\alpha\beta\gamma}^{0}=\sum_{\xi}Q_{j\xi\beta\gamma}R_{\alpha\xi}^{0}
\]
and set
\[
W_{j\alpha}(x_{2},x_{3})=\sum_{\beta,\gamma}Q_{j\alpha\beta\gamma}X_{2\beta}^{0}(x_{2})X_{3\gamma}^{0}(x_{3}).
\]
Then we define 
\[
M_{\alpha}(t,x_{1})=\sum_{\xi}X_{1\xi}(t,x_{1})R_{\xi\alpha}(t), \quad\text{ with }\quad M_{\alpha}(0,x_{1})=\sum_{\xi}X_{1\xi}^{0}(x_{1})R_{\xi\alpha}^{0},
\]
for which we obtain the one-dimensional evolution equation
\begin{equation}
\partial_t{M}_{\alpha}(t,x_{1})=-\sum_{\xi}a_{\alpha\xi}^{1}\partial_{x_1}M_{\xi}(t,x_{1})-\sum_{\xi}\left(a_{\alpha\xi}^{2}+a_{\alpha\xi}^{3}(x_{1})\right)M_{\xi}(t,x_{1})\label{eq:h6d-M}
\end{equation}
with
\[
a_{\alpha\xi}^{1}=\sum_{i,k}c_{ik}^{1;x_{1}}\la W_{i\alpha},W_{k\xi}\ra_{x_2,x_3},\qquad a_{\alpha\xi}^{2}=\sum_{i,k}\left(c_{ik}^{1;x_{2}}\la W_{i\alpha},\partial_{x_{2}}W_{k\xi}\ra_{x_2,x_3}+c_{ik}^{1;x_{3}}\la W_{i\alpha},\partial_{x_{3}}W_{k\xi}\ra_{x_2,x_3}\right)
\]
and
\[
a_{\alpha\xi}^{3}(x_{1})=\sum_{i,k}c_{ik}^{2}\cdot\la W_{i\alpha},F(M)(x_{1},\cdot,\cdot)W_{k\xi}\ra_{x_2,x_3}.
\]
This is solved up to time $\tau$. For the resulting functions $M_{\alpha}(\tau,x_{1})$ a
QR decomposition is performed to obtain orthonormal functions $X_{1\alpha}^{1}$ and the matrix $R_{\xi\alpha}^{1}$:
\[
M_{\alpha}(\tau,x_{1})=\sum_{\xi}X_{1\xi}^{1}(x_{1})R_{\xi\alpha}^{1}.
\]
Then we consider the matrix evolution equation for $R$,
\begin{equation}
\dot{R}_{\xi\alpha}(t)=\sum_{\xi',\alpha'}b_{\xi\xi'}^{1}R_{\xi'\alpha'}(t)a_{\alpha\alpha'}^{1}+\sum_{\alpha'}R_{\xi\alpha'}(t)a_{\alpha\alpha'}^{2}+\sum_{\xi',\alpha'}B_{\xi\alpha\xi'\alpha'}R_{\xi'\alpha'}(t)\label{eq:h6d-R1}
\end{equation}
with
\[
b_{\xi\xi'}^{1}=\bigl\la X_{1\xi}^{1},\partial_{x_{1}}X_{1\xi'}^{1}\bigr\ra_{x_1},\qquad B_{\xi\alpha\xi'\alpha'}=\sum_{i,k}c_{ik}^{2}\cdot\bigl\la X_{1\xi}^{1}W_{i\alpha},F(R)X_{1\xi'}^{1}W_{k\alpha'}\bigr\ra_{x}.
\]
With a low-rank representation of $F(R)$, the three-dimensional integrals in the definition of $B_{\xi\alpha\xi'\alpha'}$ can again be broken up into linear combinations of products of one-dimensional integrals.
Equation (\ref{eq:h6d-R1}) with initial value $R^{1}$ is then integrated
up to time $\tau$ to obtain $ R^{2}=R(\tau)$. Finally we obtain
\[
C_{j\alpha\beta\gamma}^{1}=\sum_{\xi}Q_{j\xi\beta\gamma} R_{\alpha\xi}^{2}.
\]

\textbf{Step 2: }We perform a QR decomposition
\[
C_{j\alpha\beta\gamma}^{1}=\sum_{\eta}Q_{j\alpha\eta\gamma}R_{\beta\eta}^{0},
\]
where $Q$ and $R$ are different from those before, but nevertheless we use the same symbols for ease of notation.
We set
\[
W_{j\alpha\beta}(x_{3})=\sum_{\gamma}Q_{j\alpha\beta\gamma}X_{3\gamma}^{0}(x_{3}).
\]
Then we define 
\[
N_{\beta}(t,x_{2})=\sum_{\eta}X_{2\eta}(t,x_{2})R_{\eta\beta}(t)
\]
for which we obtain the one-dimensional evolution equation
\begin{equation}
\partial_t{N}_{\beta}(t,x_2)=-\sum_{\eta}h_{\beta\eta}^{1}\partial_{x_{2}}N_{\eta}(t,x_{2})
-\sum_{\eta}(h_{\beta\eta}^{2}+h_{\beta\eta}^{3}(x_{2}))N_{\eta}(t,x_{2})\label{eq:h6d-N}
\end{equation}
with
\[
h_{\beta\eta}^{1}=\sum_{i,k,\alpha}c_{ik}^{1;x_{2}}\la W_{i\alpha\beta},W_{k\alpha\eta}\ra_{x_3},\qquad h_{\beta\eta}^{2}=\sum_{i,k,\alpha,\xi}c_{ik}^{1;x_{1}}b_{\alpha\xi}^{1}\la W_{i\alpha\beta},W_{k\xi\eta}\ra_{x_3}+\sum_{i,k,\alpha}c_{ik}^{1;x_{3}}\la W_{i\alpha\beta},\partial_{x_{3}}W_{k\alpha\eta}\ra_{x_3}
\]
and 
\[
h_{\beta\eta}^{3}(x_{2})=\sum_{i,k,\alpha,\xi}c_{ik}^{2}\cdot\la X_{1\alpha}^{1}W_{i\alpha\beta},F(N)(\cdot,x_{2},\cdot)X_{1\xi}^{1}W_{k\xi\eta}\ra_{x_1,x_3}.
\]
Equation (\ref{eq:h6d-N}) is then solved with initial value
\[
N_{\beta}(0,x_{2})=\sum_{\eta}X_{2\eta}^{0}(x_{2})R_{\eta\beta}^{0}
\]
up to time $\tau$. For the resulting $N_{\beta}(\tau,x_{2})$
a QR decomposition is performed to obtain $X_{2\eta}^{1}$ and $R_{\eta\beta}^{1}$.

Then we consider the matrix evolution equation for $R$,
\begin{equation}
\dot{R}_{\eta\beta}(t)=\sum_{\eta',\beta'}b_{\eta\eta'}^{2}R_{\eta'\beta'}(t)h_{\beta\beta'}^{1}+\sum_{\beta'}R_{\eta\beta'}(t)h_{\beta\beta'}^{2}+\sum_{\eta',\beta'}H_{\eta\beta\eta'\beta'}R_{\eta'\beta'}(t)
\label{eq:evol-5dh-R2}
\end{equation}
with
\[
b_{\eta\eta'}^{2}=\la X_{2\eta}^{1},\partial_{x_{2}}X_{2\eta'}^{1}\ra_{x_2},\qquad H_{\eta\beta\eta'\beta'}=\sum_{i,k,\xi,\xi'}c_{ik}^{2}\cdot\la X_{1\xi}^{1}X_{2\eta}^{1}W_{i\xi\beta},F(R)X_{1\xi'}^{1}X_{2\eta'}^{1}W_{k\xi'\beta'}\ra_x.
\]
Equation (\ref{eq:evol-5dh-R2}) with initial value $R^{1}$ is integrated
up to time $\tau$ to obtain $R^{2}$. Finally, we obtain
\[
C_{j\alpha\beta\gamma}^{2}=\sum_{\eta}Q_{j\alpha\eta\gamma}R_{\beta\eta}^{2}.
\]

\textbf{Step 3: }We perform a QR decomposition
\[
C_{j\alpha\beta\gamma}^{2}=\sum_{\zeta}Q_{j\alpha\beta\zeta}R_{\gamma\zeta}^{0}.
\]
Then we define
\[
O_{\gamma}(t,x_{3})=\sum_{\zeta}X_{3\zeta}(t,x_{3})R_{\zeta\gamma}(t)
\]
for which we obtain the evolution equation
\begin{equation}
\partial_t{O}_{\gamma}(t,x_{3})=-\sum_{\zeta}e_{\gamma\zeta}^{1}\partial_{x_{3}}O_{\zeta}(t,x_{3})-\sum_{\zeta}(e_{\gamma\zeta}^{2}+e_{\gamma\zeta}^{3}(x_{3}))O_{\zeta}(t,x_{3})\label{eq:h6d-O}
\end{equation}
with
\[
e_{\gamma\zeta}^{1}=\sum_{i,k,\alpha,\beta}c_{ik}^{1;x_{3}}Q_{i\alpha\beta\gamma}Q_{k\alpha\beta\zeta},\qquad 
e_{\gamma\zeta}^{2}=\sum_{i,k,\alpha,\xi,\beta}c_{ik}^{1;x_{1}}b_{\alpha\xi}^{1}Q_{i\alpha\beta\gamma}Q_{k\xi\beta\zeta}+\sum_{i,k,\alpha,\beta,\eta}c_{ik}^{1;x_{2}}b_{\beta\eta}^{2}Q_{i\alpha\beta\gamma}Q_{k\alpha\eta\zeta}
\]
and
\[
e_{\gamma\zeta}^{3}(x_{3})=\sum_{i,k,\alpha,\xi,\beta,\eta}c_{ik}^{2}\cdot\la X_{1\alpha}X_{2\beta},F(O)(\cdot,\cdot,x_{3})X_{1\xi}X_{2\eta}\ra_{x_1,x_2}Q_{i\alpha\beta\gamma}Q_{k\xi\eta\zeta}.
\]
Equation (\ref{eq:h6d-O}) is then solved with initial value 
\[
O_{\gamma}(0,x_{3})=\sum_{\zeta}X_{3\zeta}^{0}(x_{3})R_{\zeta\gamma}^{0}
\]
up to time $\tau$. For the resulting $O_{\gamma}(\tau,x_{2})$ a QR
decomposition is performed to obtain $X_{3\zeta}^{1}$ and $R_{\zeta\gamma}^{1}$.

Then we consider the evolution equation for $R$,
\begin{equation}
\dot{R}_{\zeta\gamma}(t)=\sum_{\zeta',\gamma'}b_{\zeta\zeta'}^{3}R_{\zeta'\gamma'}(t)e_{\gamma\gamma'}^{1}+\sum_{\gamma'}R_{\zeta\gamma'}(t)e_{\gamma\gamma'}^{2}+\sum_{\zeta',\gamma'}G_{\zeta\gamma\zeta'\gamma'}R_{\zeta'\gamma'}(t)\label{eq:h6d-R3}
\end{equation}
with
\begin{align*}
b_{\gamma\zeta'}^{3} & =(X_{3\gamma},\partial_{x_{3}}X_{3\zeta'}),\qquad G_{\zeta\gamma\zeta'\gamma'}=\sum_{i,k,\alpha,\beta,\alpha',\beta'}c_{ik}^{2}\cdot\bigl\la X_{1\alpha}^{1}X_{2\beta}^{1}X_{3\zeta}^{1},F(R)X_{1\alpha'}^{1}X_{2\beta'}^{1}X_{3\zeta'}^{1}\bigr\ra_x
Q_{i\alpha\beta\gamma}Q_{k\alpha'\beta'\gamma'}.
\end{align*}
Equation (\ref{eq:h6d-R3}) with initial value $R^{1}$ is integrated
up to time $\tau$ to obtain $R^{2}$. Finally we obtain
\[
C_{j\alpha\beta\gamma}^{3}=\sum_{\gamma}Q_{j\alpha\beta\zeta}R_{\gamma\zeta}^{2}.
\]

\textbf{Step 4: }The final step directly updates $C$. The corresponding
evolution equation is
\begin{equation}
\dot{C}_{i\alpha\beta\gamma}(t)=\sum_{k,\xi}c_{ik}^{1;x_{1}}b_{\alpha\xi}^{1}C_{k\xi\beta\gamma}(t)+\sum_{k,\eta}c_{ik}^{1;x_{2}}b_{\beta\eta}^{2}C_{k\alpha\eta\gamma}(t)+\sum_{k,\zeta}c_{ik}^{1;x_{3}}b_{\gamma\zeta}^{3}C_{k\alpha\beta\zeta}(t)+\sum_{k,\xi,\eta,\zeta}(c_{ik}^{2}\cdot e_{\alpha\beta\gamma\xi\eta\zeta})C_{k\xi\eta\zeta}(t)\label{eq:h6d-evol-C}
\end{equation}
with
\[
e_{\alpha\beta\gamma\xi\eta\zeta}=\bigl\la X_{1\alpha}^{1}X_{2\beta}^{1}X_{3\gamma}^{1},F(R)X_{1\xi}^{1}X_{2\eta}^{1}X_{3\zeta}^{1}\bigr\ra_x.
\]
Equation (\ref{eq:h6d-evol-C}) with initial value $C^{3}$ is then
solved up to time $\tau$ to obtain $C^{4}$. This completes a full time step of algorithm.
The output is 
\[
K_{j}(\tau,x)\approx K_{j}^{1}(x)=\sum_{\alpha,\beta,\gamma}C_{j\alpha\beta\gamma}^{4}X_{1\alpha}^{1}(x_{1})X_{2\beta}^{1}(x_{2})X_{3\gamma}^{1}(x_{3}).
\]
At the end, a factorization into an orthogonalized tensor $C$ and a matrix $S$ is obtained as in Algorithm~2, simply by adding a subscript $\gamma$ to every appearance of $\alpha,\beta$ in that algorithm.

\section{Numerical results\label{sec:Numerical-results}}

\change{In this section we show numerical results for the low-rank and the hierarchical low-rank splitting algorithm. In both cases the spectral method described in section \ref{subsec:FFT} is used. Numerical results for linear Landau damping, for a two-stream instability and a plasma echo will be presented.}

\subsection{Linear Landau damping\label{subsec:ll}}

\begin{figure}
\centering{}\includegraphics[width=16cm]{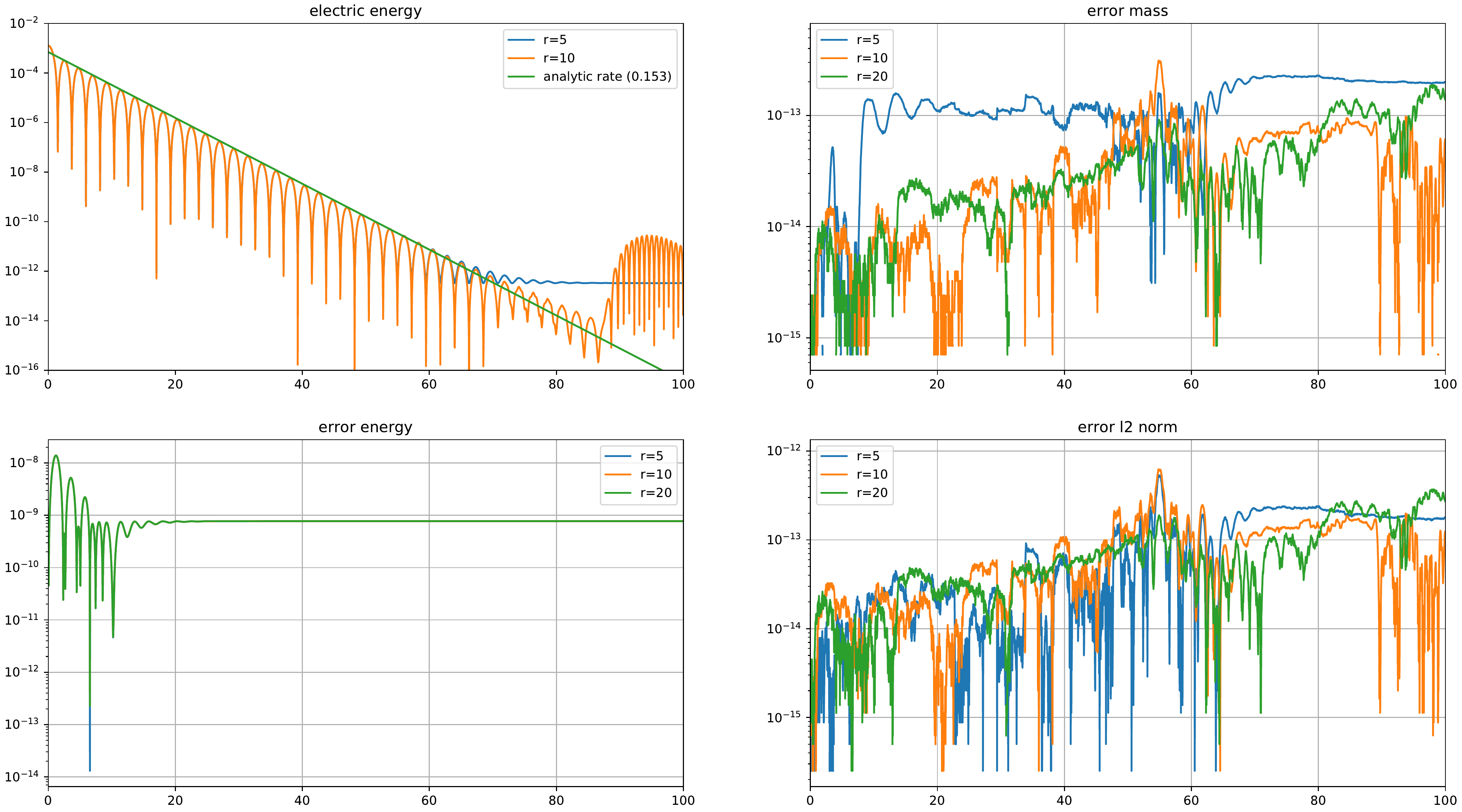}\caption{Numerical simulations of linear Landau damping with the low-rank/FFT
algorithm and various ranks are shown (note that the plots for $r=10$
and $r=20$ in the top left figure are indistinguishable). The Strang
splitting algorithm with a time step size $\tau=0.025$ is employed.
In the $x$-direction $64$ grid points are used, while in the $v$-direction
$256$ grid points are used. \label{fig:ll2d}}
\end{figure}

First, let us consider the classic two-dimensional Landau damping.
We set $\Omega=\Omega_x\times\Omega_v=(0,4\pi)\times(-6,6)$ and impose the initial
value
\[
f_{0}(x,v)=\frac{1}{\sqrt{2\pi}}\mathrm{e}^{-v^{2}/2}(1+\alpha\cos(kx)),
\]
where we have chosen $\alpha=10^{-2}$ and $k=\tfrac{1}{2}$. In both
the $x$- and $v$-direction periodic boundary conditions are imposed.
It can be shown by a linear analysis that the decay rate of the electric
energy is given by $\gamma\approx-0.153$. The so obtained decay rate
has been verified by a host of numerical simulations presented
in the literature. The numerical results obtained using the low-rank/FFT
algorithm proposed in Section \ref{subsec:FFT} are shown in Figure
\ref{fig:ll2d}.

 We observe that choosing the rank $r=5$ is sufficient to obtain a numerical
solution which very closely matches the analytic result. Although
we see that the numerical method does not conserve energy up to machine
precision, the error is very small (on the order of $10^{-8}$). In
addition, the errors in mass and in the $L^{2}$ norm are indistinguishable
from machine precision. 

\begin{figure}
\centering{}\includegraphics[width=16cm]{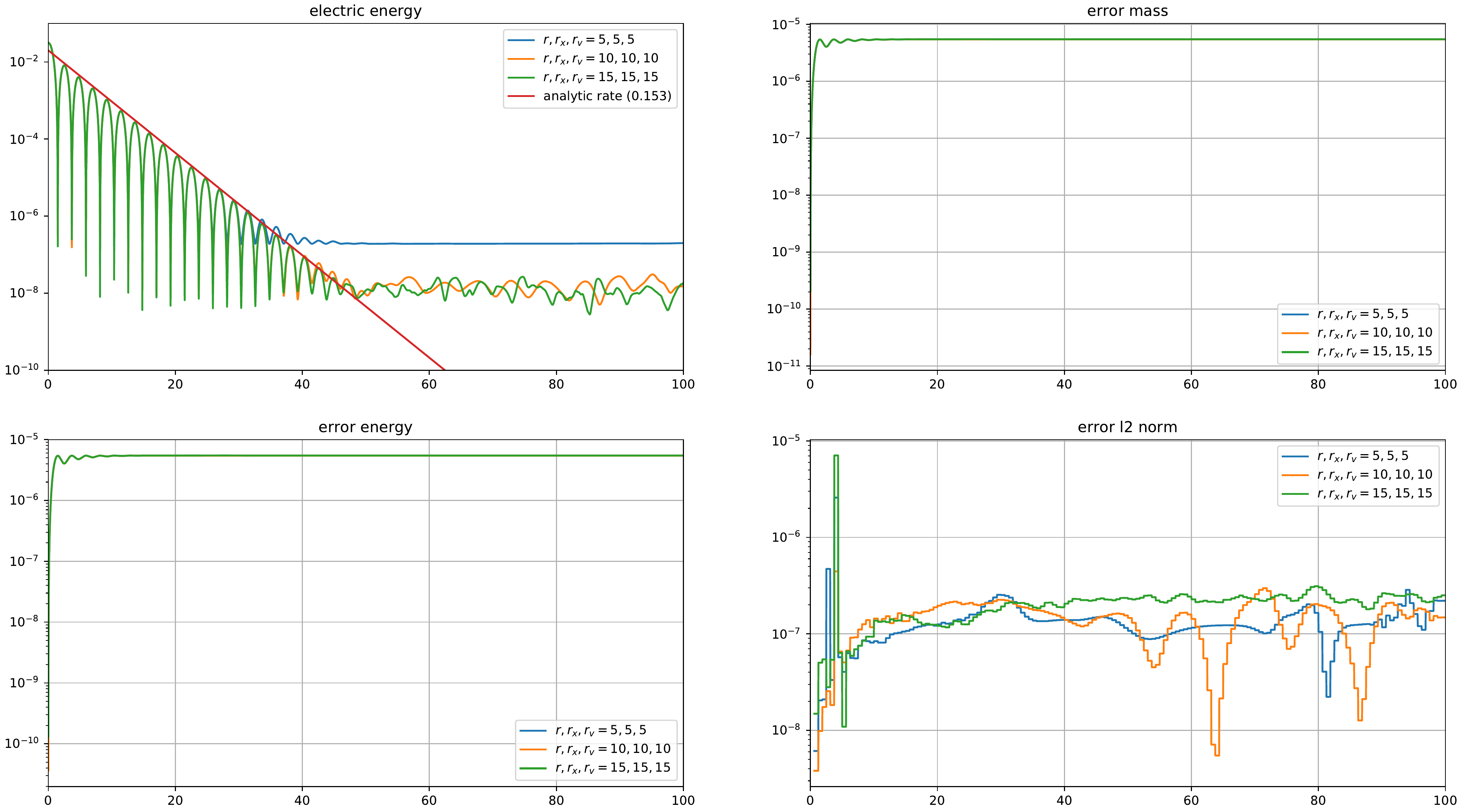}\caption{Numerical simulations of four-dimensional linear Landau damping with
the hierarchical low-rank/FFT algorithm and various ranks are shown
(note that the error in mass and energy is almost identical for all
three configurations). The Lie--Trotter splitting algorithm with a time step
size $\tau=0.00625$ is employed. In the $x$-direction $64$ grid
points are used, while in the $v$-direction $256$ grid points are
used. The $L^{2}$ norm is only
computed for every $40$th time step. \label{fig:ll4d}}
\end{figure}

Next, we turn our attention to a four-dimensional problem. That is,
we set $\Omega=(0,4\pi)^{2}\times(-6,6)^{2}$ and impose the initial
value

\[
f_{0}(x,y,v,w)=\frac{1}{2\pi}\mathrm{e}^{-(v^{2}+w^{2})/2}(1+\alpha\cos(k_{1}x)+\alpha\cos(k_{2}y)),
\]
where we have chosen $\alpha=10^{-2}$ and $k_{1}=k_{2}=\tfrac{1}{2}$.
As before, periodic boundary conditions are chosen for all directions.
Although the problem is now set in four dimensions, the behavior of
the electric energy is almost identical to the two-dimensional problem.
The numerical results obtained using the hierarchical low-rank algorithm
proposed in section \ref{sec:The-hierarchical-algorithm} are shown
in Figure \ref{fig:ll4d}. 
 We observe that the simulation with rank $(r,r_{x},r_{v})=(5,5,5)$
is already able to correctly predict the decay rate. Considering the
simulation with $(r,r_{x},r_{v})=10$ shows a reduction in the electric
energy to $10^{-8}$ which is approximately two orders of magnitude
better compared to the configuration with rank $(r,r_{x},r_{v})=(5,5,5)$.
In all configurations mass, energy, and the $L^{2}$ norm are conserved
up to an error less than $10^{-6}$. 

\begin{figure}
    \centering{}\includegraphics[width=16cm]{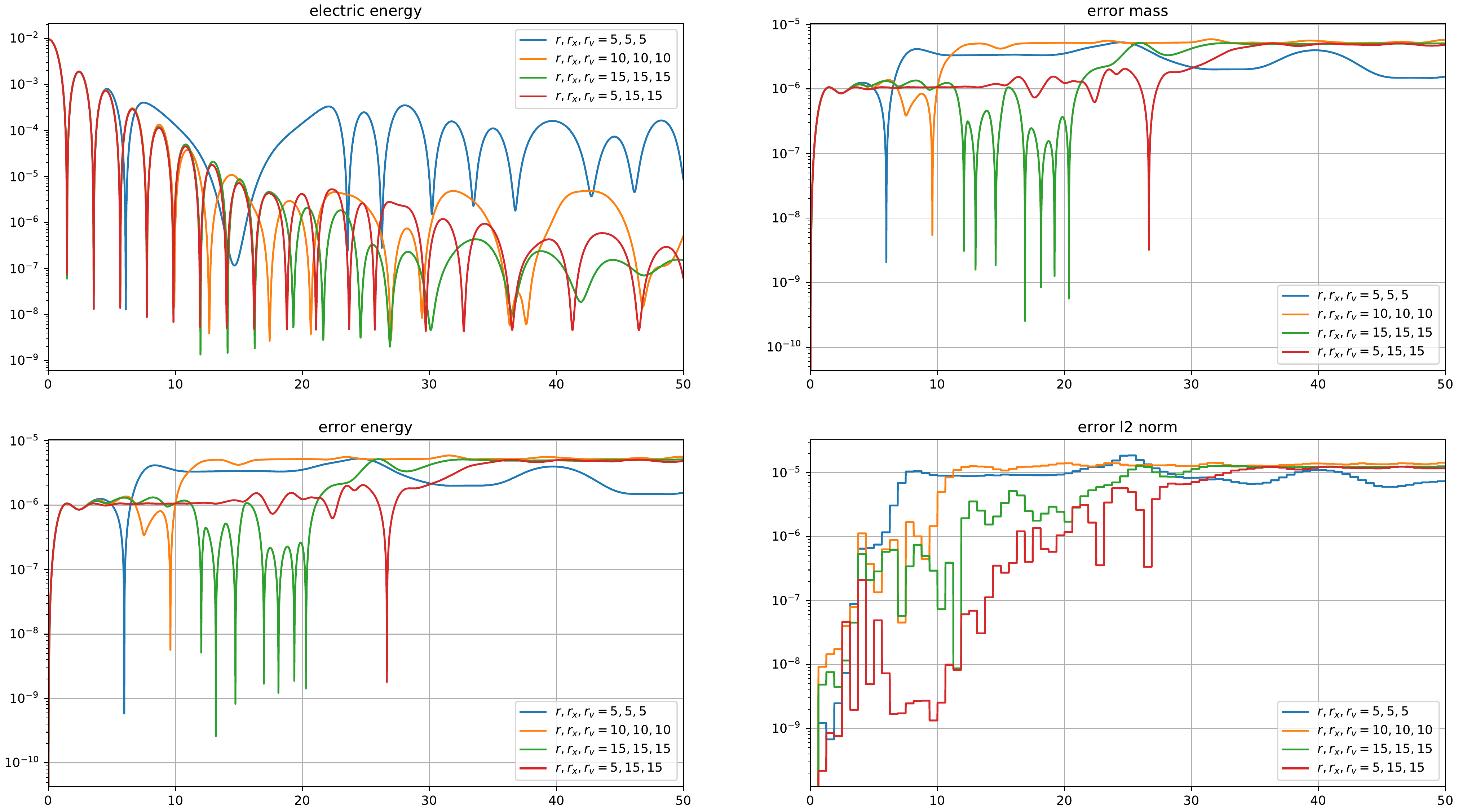}\caption{\change{Numerical simulations of the four-dimensional non-aligned linear Landau damping with the hierarchical low-rank/FFT algorithm and various ranks are shown
(note that the error in mass and energy is almost identical for all
three configurations). The Lie--Trotter splitting algorithm with a time step
size $\tau=0.00625$ is employed. In the $x$-direction $64$ grid
points are used, while in the $v$-direction $256$ grid points are
    used. The $L^{2}$ norm is only computed for every $100$th time step.} \label{fig:ll4d_notaligned}}
\end{figure}

\change{To conclude the discussion on Landau damping, we consider a setting in which the perturbation is not aligned to the coordinate axis. More specifically, we consider $\Omega=(0,5\pi)^2\times (-6,6)^2$ and impose the initial value 
\[ f_{0}(x,y,v,w)=\frac{1}{2\pi}\mathrm{e}^{-(v^{2}+w^{2})/2}(1+\alpha\cos(k_{1}x)\cos(k_{2}y)), \]
where we have chosen $\alpha=10^{-2}$ and $k_{1}=k_{2}=0.4$. As before, periodic boundary conditions are chosen for all directions. The numerical results obtained using the hierarchical low-rank algorithm proposed in section \ref{sec:The-hierarchical-algorithm} are shown in Figure \ref{fig:ll4d_notaligned}. We observe that the simulation with rank $(r,r_x,r_v)=(10,10,10)$ and $(r,r_x,r_v)=(15,15,15)$ show a reduction of the electric energy to $5\cdot 10^{-6}$ and $5\cdot 10^{-7}$, respectively. We also note that $r$ can be chosen significantly smaller than $r_x$ and $r_v$, while still obtaining good results. For all configurations the error in mass, energy, and $L^2$ norm is below $10^{-5}$.}

\subsection{Two-stream instability}

Our second numerical example is the so-called two-stream instability.
Here we have two beams propagating in opposite directions. This
setup is an unstable equilibrium and small perturbations in the initial
particle-density function eventually force the electric energy to
increase exponentially. This is called the linear regime. At some
later time saturation sets in (the nonlinear regime). This phase is
characterized by nearly constant electric energy and significant filamentation
of the phase space.

\begin{figure}
\centering{}\includegraphics[width=16cm]{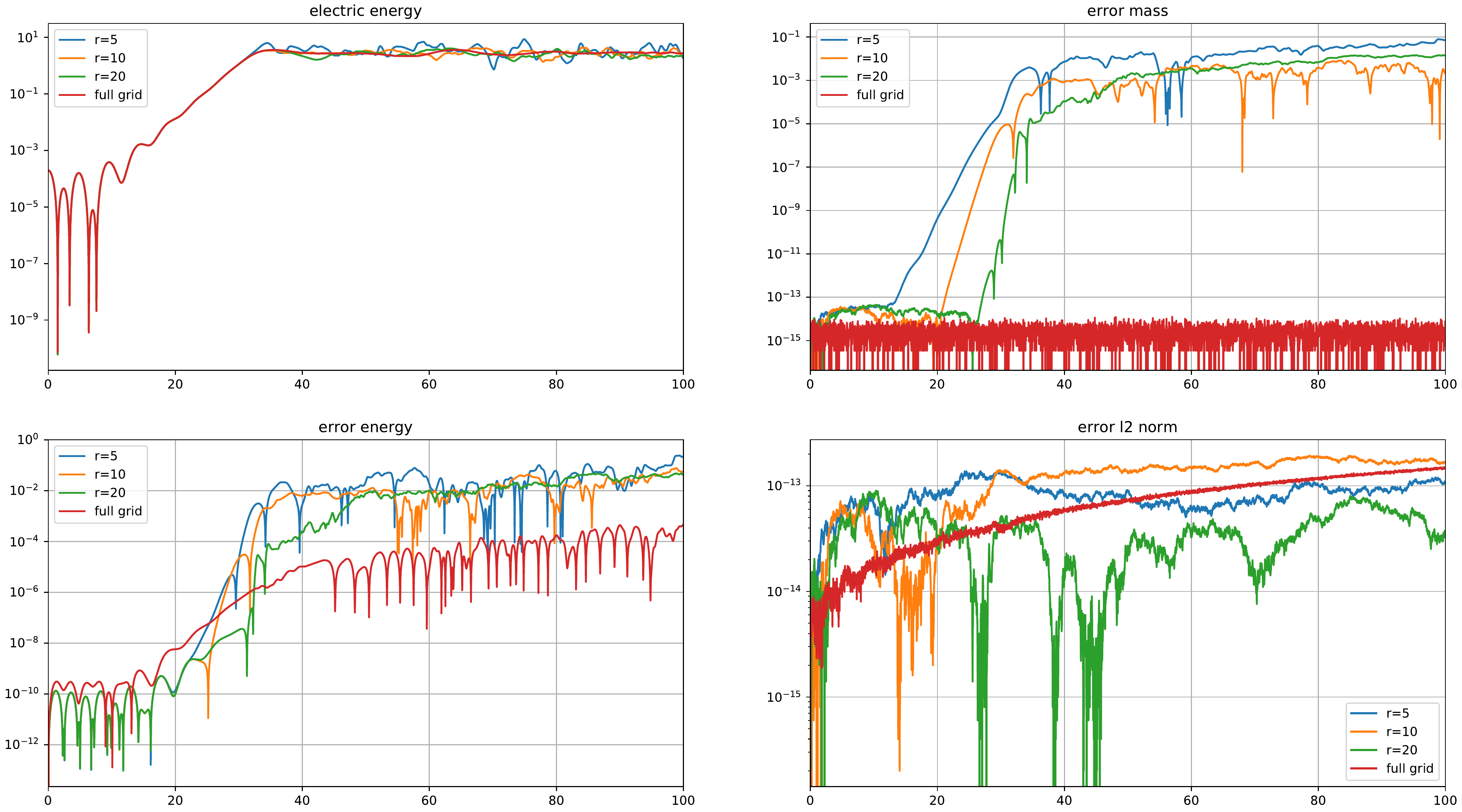}\caption{Numerical simulations of the two-dimensional two-stream instability
with the low-rank/FFT algorithm and various ranks are shown. The Strang
splitting algorithm with a time step size $\tau=0.025$ is employed.
In both the $x$ and $v$-directions $128$ grid points are used.
As comparison a direct Eulerian simulation with a splitting/discontinuous
Galerkin semi-Lagrangian scheme is also shown. \label{fig:ts2d}}
\end{figure}

First we consider the two-dimensional case. We thus set \change{$\Omega=(0,10\pi)\times(-9,9)$} and impose the initial condition
\[
f_{0}(x,v)=\frac{1}{2\sqrt{2\pi}}\left(\mathrm{e}^{-(v-v_{0})^{2}/2}+\mathrm{e}^{-(v+v_{0})^{2}/2}\right)(1+\alpha\cos(kx)),
\]
where $\alpha=10^{-3}$, $k=\tfrac{1}{5}$, and $v_{0}=2.4$. As before,
periodic boundary conditions are used. The corresponding numerical
simulations obtained using the low-rank/FFT algorithm proposed in
Section \ref{subsec:FFT} are shown in Figure \ref{fig:ts2d}. 
 In the linear regime excellent agreement between the direct Eulerian
simulation and the solution with rank $r=5$ is observed. As we enter
the nonlinear regime all the solutions become somewhat distinct. However,
due to the chaotic nature of this regime this is not surprising. The
figure of merit we are looking at in the nonlinear regime is if the
numerical scheme is able to keep the electric energy approximately
constant. All configurations starting from rank $r=10$ do this very well.

To further evaluate the quality of the numerical solution, let us
consider the physical invariants. By looking at Figure \ref{fig:ts2d},
it becomes clear that from this standpoint the two-stream instability
is a significantly more challenging problem than the linear
Landau damping considered in the previous section. However, it should
be noted that most of the error in mass and energy is accrued
during the exponential growth of the electric energy but remains essentially
unchanged in the nonlinear phase. In particular, there is no observable
long-term drift in the error of any of the invariants. With respect
to energy, the low-rank approximation with $r=10$ is 
worse by approximately two orders of magnitude compared to the direct Eulerian simulation
(which has been conducted using a FFT based implementation).
The $L^{2}$ norm is conserved up to machine precision.

\begin{figure}
\centering{}\includegraphics[width=16cm]{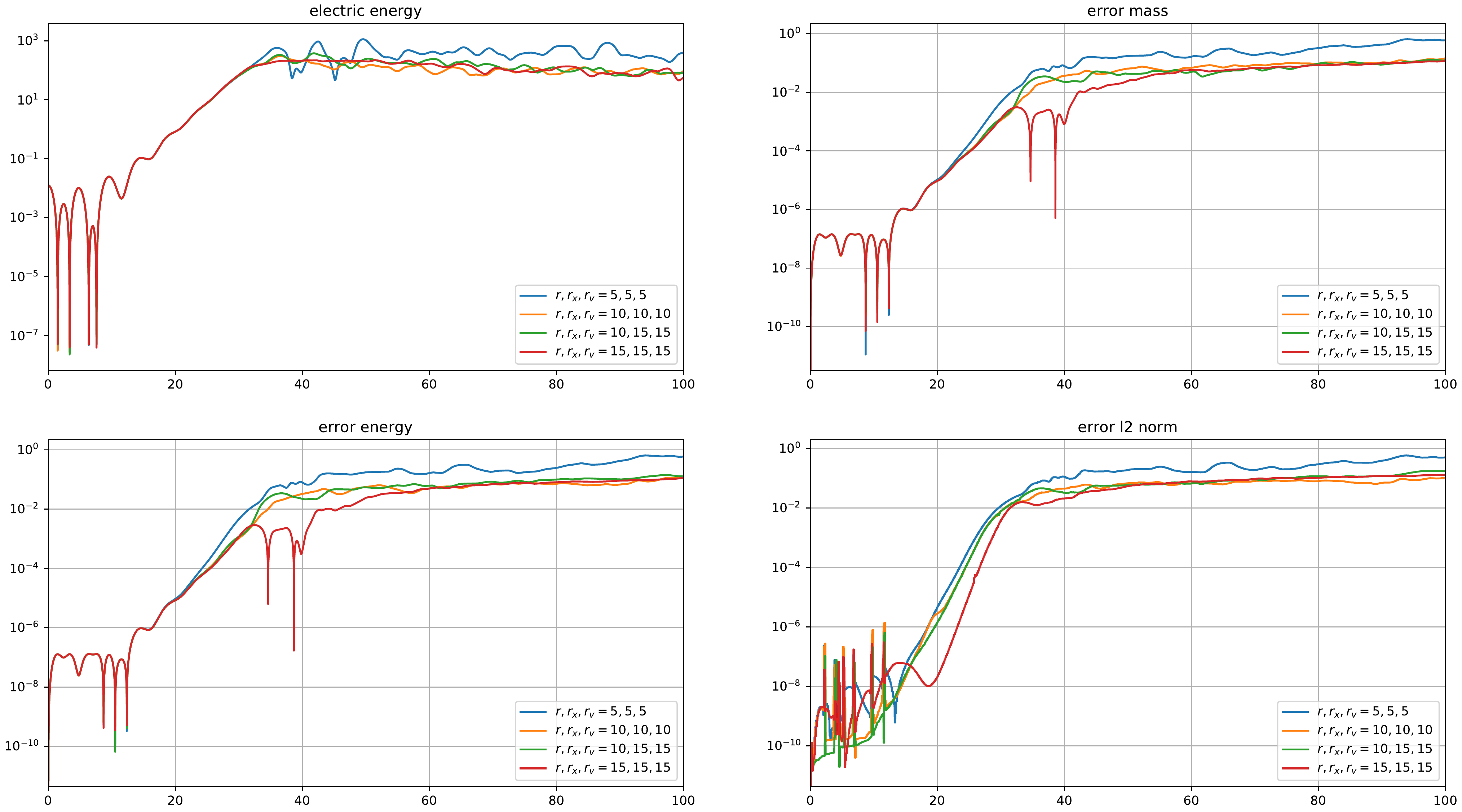}\caption{Numerical simulations of the four-dimensional two-stream instability
with the hierarchical low-rank/FFT algorithm and various ranks are
shown. The Lie--Trotter projector-splitting algorithm with a time step size $\tau=0.00625$
is employed. In both the $x$- and $v$-directions $128$ grid points
are used. The $L^{2}$ norm is only
computed for every $100$th time step. \label{fig:ts4d}}
\end{figure}

Second, we turn our attention to the four-dimensional case. We set
\change{$\Omega=(0,10\pi)^{2}\times(-9,9)^{2}$} and impose the initial
condition
\[
f_{0}(x,y,v,w)=\frac{1}{8\pi}\left(\mathrm{e}^{-(v-v_{0})^{2}/2}+\mathrm{e}^{-(v+v_{0})^{2}/2}\right)\left(\mathrm{e}^{-(v-w_{0})^{2}/2}+\mathrm{e}^{-(v+w_{0})^{2}/2}\right)(1+\alpha_{1}\cos(k_{1}x)+\alpha_{2}\cos(k_{2}y))
\]
with $\alpha_{1}=\alpha_{2}=10^{-3}$, $k_{1}=k_{2}=\tfrac{1}{5}$,
$v_{0}=w_{0}=2.4$. As before, periodic boundary conditions are employed
in all directions. We note that this configuration corresponds to
four propagating beams. The numerical results obtained using the hierarchical
low-rank algorithm proposed in section \ref{sec:The-hierarchical-algorithm}
are shown in Figure \ref{fig:ts4d}.
 The obtained results are similar to the two-dimensional case. In
particular, in the linear regime we observe excellent agreement even
for $(r,r_{x},r_{v})=(5,5,5)$. The errors in mass, energy, and $L^{2}$
norm are approximately one order of magnitude better for the configuration
with $(r,r_{x},r_{v})=(10,10,10)$ compared to $(r,r_{x},r_{v})=(5,5,5)$.

\subsection{Plasma echo}

\change{There are problems that need to resolve Landau damping for relatively long times.} It is well known that due to the recurrence effect the number of grid points $n_v$ in  each coordinate direction of $v$ has to scale as $n_{v}\propto T$, where $T$ is the final time of the simulation (although this can be somewhat alleviated by using filamentation filtration techniques; see, for example, \cite{klimas1994,einkemmer2014}).
\change{Due to the high number of grid points required to conduct such simulation the proposed algorithm can be extremely efficient in this context. We will now consider such an example, the plasma echo phenomenon (see, for example, \cite{gould1967,hou2011,einkemmer2014}).}

\change{We consider the domain $(0,100)\times(-8,8)$ and impose the initial condition
\[ f_0(x,v)=\frac{1}{2\pi} \mathrm{e}^{-v^2/2} ( 1+ \alpha \cos (k_1 x) ), \]
where $\alpha=10^{-3}$ and $k_1=12\pi/100$. As before, periodic boundary conditions are used. The initial perturbation in the electric field is damped away by Landau damping. Then at time $t_2=200$ we add another perturbation of the form\[ \frac{\alpha}{2\pi} \mathrm{e}^{-v^2/2}  \cos (k_2 x) , \]
where $k_2 = 24\pi/100$. The corresponding perturbation in the electric field is, again, damped away. However, the particle-density function retains all the information regarding the two perturbations. As a consequence, at later times an echo (a peak in the amplitude of the electric field) occurs. The time of the echo depends on the wave numbers of the two perturbations. In the present case we expect a primary echo at $t=400$ and a secondary echo at $t=800$. This setup has been considered in \cite{hou2011} and \cite{einkemmer2014}.}

\begin{figure}
    \centering{}\includegraphics[width=16cm]{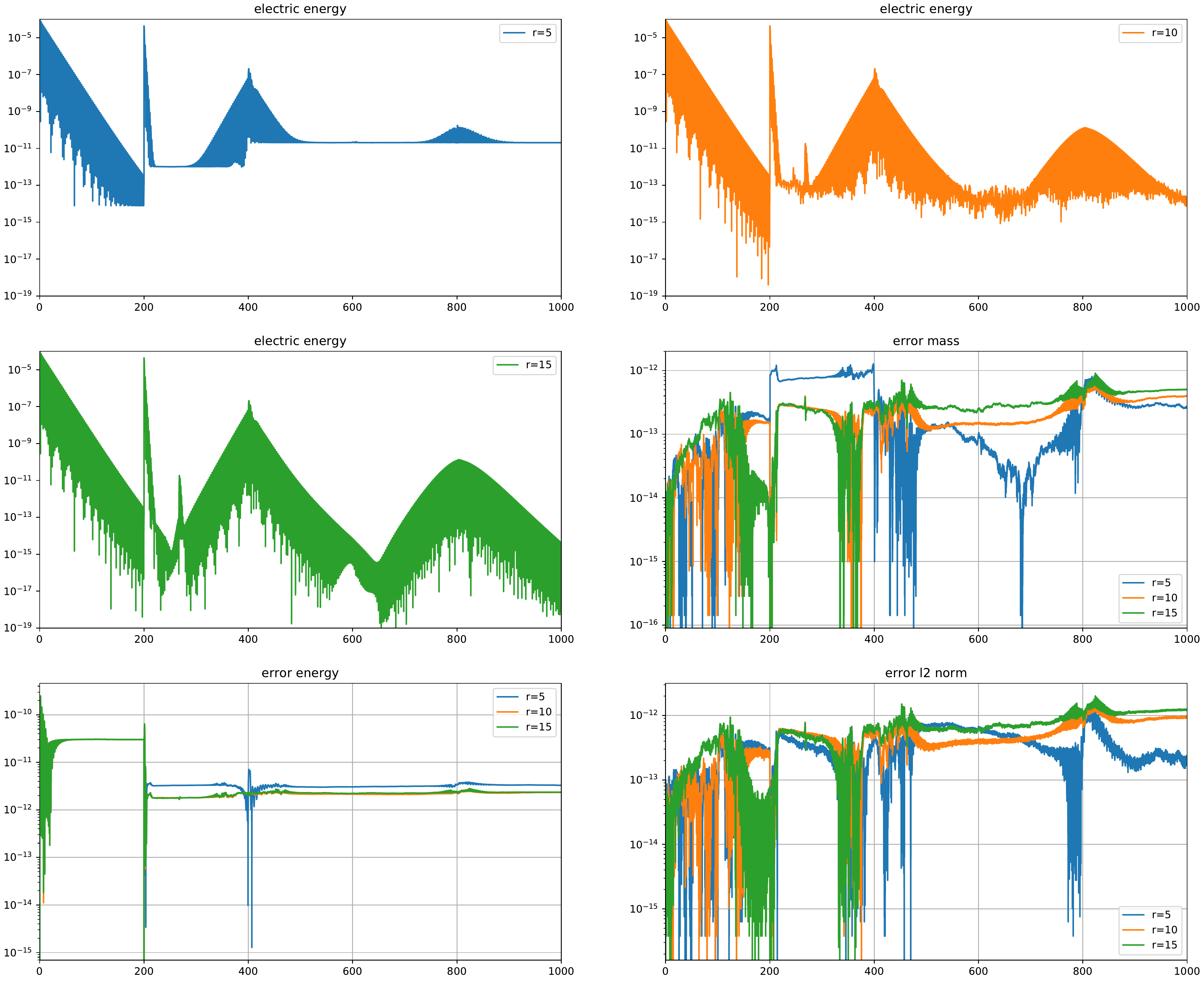}\caption{\change{Numerical simulations of the two-dimensional plasma echo with the low-rank/FFT algorithm and various ranks are shown. The Strang projector-splitting algorithm with a time step size $\tau=0.025$ is employed. In the $x$- and $v$-direction $512$ and $4096$ grid points are used, respectively. Note that due to the second perturbation, we restart the computation of the errors in mass, energy, and $L^2$ norm at $t=200$.} \label{fig:echo}}
\end{figure}

\change{The corresponding numerical simulations are shown in Figure \ref{fig:echo}. Even for $r=5$ both the primary and the secondary echo are clearly resolved. This holds true even though the electric energy is only resolved up to an amplitude of approximately $10^{-11}$ and the magnitude of the secondary echo is only slightly above that threshold. Increasing the rank (to $r=10$) then allows us to resolve the growth and decay of the electric energy associated with the secondary echo as well. Since this problem requires only a small rank but a large number of grid points ($4096$ grid points are used in the velocity direction and $512$ grid points are used in the space direction), the computational effort required is reduced by a significant margin. Let us also note that, for all configurations considered here, the conservation of mass, energy, and $L^2$ norm is excellent (below $10^{-10}$).}

\change{Overall, the results presented here for Landau damping, the two-stream instability, and the plasma echo, show that significant reduction in computational effort can be obtained compared with performing a direct simulation. In fact, all simulations have been
performed on a somewhat outdated workstation and no effort has been made to parallelize the code.}

\bigskip\bigskip
\pagebreak[3]
\bibliographystyle{plain}
\bibliography{vlasov-general,vlasov-lowrank}

\end{document}